\documentclass[11pt]{article}

\usepackage{color}
\usepackage{latexsym}
\usepackage{dsfont}
\usepackage{amssymb}
\usepackage{graphicx}
\usepackage{amsmath, amsfonts,amssymb,theorem,euscript,array,enumerate,amsfonts,mathrsfs}
\usepackage{hyperref}

\newtheorem{Theorem}{Theorem}[part]
\newtheorem{Definition}{Definition}[part]
\newtheorem{Proposition}{Proposition}[part]

\newtheorem{Lemma}{Lemma}[part]

\newtheorem{Remark}{Remark}[part]

\def \trans{^{\scriptscriptstyle{\intercal}}}

\def \Frac{\displaystyle\frac}
\def \Inf{\displaystyle\inf}
\def \Sup{\displaystyle\sup}

\def \b1{\bf{1}}
\def \bA{\bf{A}}

\def \N{\mathbb{N}}
\def \R{\mathbb{R}}
\def \L{\mathbb{L}}
\def \M{\mathbb{M}}

\def \E{\mathbb{E}}
\def \F{\mathbb{F}}

\def \P{\mathbb{P}}
\def \S{\mathbb{S}}
\def \Q{\mathbb{Q}}

\def\esssup_#1{\underset{#1}{\mathrm{ess\,sup\, }}}
\def\essinf_#1{\underset{#1}{\mathrm{ess\,inf\, }}}

\def\argmin_#1{\underset{#1}{\mathrm{argmin\, }}}

\def \Ac{{\cal A}}
\def \Bc{{\cal B}}
\def \Cc{{\cal C}}

\def \Ec{{\cal E}}
\def \Fc{{\cal F}}
\def \Gc{{\cal G}}
\def \Hc{{\cal H}}

\def \Kc{{\cal K}}
\def \Lc{{\cal L}}
\def \Pc{{\cal P}}

\def \Tc{{\cal T}}

\def \Wc{{\cal W}}

\def \eps{\varepsilon}

\def \ep{\hbox{ }\hfill$\Box$}

\def\Dt#1{\Frac{\partial #1}{\partial t}}

\def\reff#1{{\rm(\ref{#1})}}

\def\beqs{\begin{eqnarray*}}
\def\enqs{\end{eqnarray*}}
\def\beq{\begin{eqnarray}}
\def\enq{\end{eqnarray}}

\begin{document}

 \title{Dynamic programming for optimal control of stochastic McKean-Vlasov dynamics 
 %and applications
 \thanks{This work is part of the ANR project CAESARS (ANR-15-CE05-0024). 
 We would like to thank the referees for their suggestions which help us to improve the paper.
 }
 %\thanks{We would like to thank E. Bayraktar, J.F. Chasagneux and A. Cosso for discussions during the preparation of this paper.}
 }

\author{Huy\^en PHAM
\\\small  Laboratoire de Probabilit\'es et
 \\\small  Mod\`eles Al\'eatoires, CNRS, UMR 7599
 \\\small  Universit\'e Paris Diderot
 \\\small  pham at math.univ-paris-diderot.fr
\\\small  and CREST-ENSAE
\and
Xiaoli WEI
\\\small  Laboratoire de Probabilit\'es et
 \\\small  Mod\`eles Al\'eatoires, CNRS, UMR 7599
 \\\small  Universit\'e Paris Diderot
 \\\small  tyswxl at gmail.com
}

%\author{Huy{\^e}n PHAM\thanks{Laboratoire de Probabilit\'es et Mod\`eles Al\'eatoires, CNRS, UMR 7599, Universit{\'e} Paris Diderot, and
 %              CREST-ENSAE,  \sf pham at math.univ-paris-diderot.fr} ~~~
%		Xiaoli WEI\thanks{Laboratoire de Probabilit\'es et Mod\`eles Al\'eatoires, CNRS, UMR 7599, Universit\'e Paris Diderot,
%		\sf tyswxl at gmail.com}
 %            }

\maketitle

\date{}

\begin{abstract}
We study the optimal control of general stochastic McKean-Vlasov equation. Such problem is motivated originally from the asymptotic formulation of cooperative equilibrium 
for a large population of particles (players) in mean-field interaction under common noise. 
%It also arises from partial observation control problem after a change of probability reference. 
Our first main result is to state a dynamic programming principle for the value function in the Wasserstein space of probability measures, which is proved from   
a flow property of the conditional law of the controlled state process.
Next, by relying on the notion of differentiability with respect to  probability measures due to P.L. Lions \cite{lio12}, and It\^o's formula along a flow of 
conditional measures, we derive the dynamic programming Hamilton-Jacobi-Bellman equation, and prove the viscosity  property together with a uniqueness result for the value function. 
Finally, we solve explicitly the linear-quadratic stochastic McKean-Vlasov control problem and give an application to an interbank systemic risk model with common noise.  
\end{abstract}

\vspace{5mm}

 %60K35, 49L20
%60G40, 91A05,  49L20,  49L25.

\vspace{5mm}


\newpage

 \section{Introduction}

\setcounter{equation}{0} \setcounter{Assumption}{0}
\setcounter{Theorem}{0} \setcounter{Proposition}{0}
\setcounter{Corollary}{0} \setcounter{Lemma}{0}
\setcounter{Definition}{0} \setcounter{Remark}{0}

Let us consider the controlled McKean-Vlasov dynamics in $\R^d$  given by
\beq \label{stoMcKean}
dX_t &=& b(X_t,\P_{X_t}^{W^0},\alpha_t) dt + \sigma(X_t,\P_{X_t}^{W^0},\alpha_t) dB_t + \sigma_0(X_t,\P_{X_t}^{W^0},\alpha_t) dW^0_t,
\enq
where $B,W^0$ are two independent Brownian motions on some complete probability space $(\Omega,\Fc,\P)$, $\P_{X_t}^{W^0}$ denotes the conditional distribution of $X_t$ given $W^0$ (or equivalently given $\Fc_t^0$ where $\F^0$ $=$ $(\Fc_t^0)_{t\geq 0}$ is the natural filtration generated by $W^0$), 
valued in $\Pc(\R^d)$ the set of probability measures on $\R^d$, 
and the control $\alpha$ is an $\F^0$-progressive process valued in some Polish space $\bA$. When there is no control, 
the dynamics \reff{stoMcKean} is sometimes called stochastic McKean-Vlasov equation (see \cite{dawvai95}), where the term ``stochastic"  refers to the presence of the random noise caused by the Brownian motion $W^0$ w.r.t.  a McKean-Vlasov equation when $\sigma_0$ $=$ $0$, and for which coefficients depend on the (deterministic) marginal distribution $\P_{X_t}$.  One also uses the terminology  conditional mean-field stochastic differential equation (CMFSDE) to emphasize the dependence of the coefficients on the conditional law with respect to the random noise, and such CMFSDE was studied in \cite{carzhu14}, and more generally in \cite{buclima15}. 
In this context, the control problem is to minimize over $\alpha$ a cost functional of the form:
\beq \label{costMcKean}
J(\alpha) &=& \E \Big[ \int_0^T f(X_t,\P_{X_t}^{W^0},\alpha_t) dt + g(X_T,\P_{X_T}^{W^0}) \Big]. 
\enq

The motivation and applications for the study of such stochastic control problem, referred to alternatively  as control of stochastic McKean-Vlasov dynamics, or 
stochastic control of conditional McKean-Vlasov equation,  comes mainly from   the  {\it  McKean-Vlasov control problem with common noise}, that we briefly describe now: 
we consider a system of controlled individuals (referred also  to as particles or players) in mutual interaction, 
where the dynamics of the state process $X^i$  
of player $i$ $\in$ $\{1,\ldots,N\}$ is governed by 
\beqs
dX_t^i &=& \tilde b(X_t^i, \bar\rho_t^N,\tilde \alpha_t^i) dt 
+  \tilde\sigma(X_t^i, \bar\rho_t^N,\tilde\alpha_t^i) dB_t^i +  \tilde\sigma_0(X_t^i, \bar\rho_t^N,\tilde\alpha_t^i) dW_t^0.   
\enqs
Here, the Wiener process $W^0$  accounts for the common random environment in which all the individuals evolve, called common noise, and 
$B^1,\ldots,B^N$ are independent Brownian motions, independent of $W^0$,  called idiosyncratic noises.  The particles are in interaction of mean-field type in the sense that any any time $t$, the coefficients $\tilde b$, $\tilde\sigma$, $\tilde\sigma_0$ of their state process 
depend on the empirical  distribution of all individual states
\beqs
\bar\rho_t^N &=& \frac{1}{N} \sum_{i=1}^N \delta_{X_t^i}. 
\enqs
The processes $(\tilde\alpha_t^i)_{t\geq 0}$, $i$ $=$ $1,\ldots,N$,  are  in general progressively measurable w.r.t. the filtration generated by 
$B^1,\ldots,B^N,W^0$,  valued in some subset $A$  of a Euclidian space,  and represent the control processes  of the players with cost functionals: 
\beqs
J^i(\tilde\alpha^1,\ldots,\tilde\alpha^n) &=& \E \Big[ \int_0^T \tilde f(X_t^i,\bar\rho_t^N,\tilde\alpha_t^i) dt + g(X_T^i,\bar\rho_T^N) \Big]. 
\enqs 
For this $N$-player stochastic differential game, one looks for equilibriums, and different notions  may be considered. Classically,  
the search for a consensus among the players  leads to the concept of  Nash equilibrium where  each player minimizes its own cost  functional, and the goal is to find a $N$-tuple control strategy for which there is no interest 
for any player to leave from this consensus state. The asymptotic formulation of this Nash equili\-brium when the number of players $N$ goes  to infinity leads to the (now well-known) theory of mean-field games (MFG) pioneered in the 
works by Lasry and Lions \cite{laslio07}, and Huang, Malham\'e and Caines \cite{huaetal06}.  In this framework,  the analysis is reduced to the problem of a single representative player in interaction with the theoretical distribution of the whole population by the propagation of chaos phenomenon, who first solves a control problem by freezing a probability law  in the coefficients of her/his state process and cost function, and then has to find a fixed point probability measure that matches the distribution of her/his optimal state process. The case of MFG with common noise has been recently studied in   \cite{ahu14} and  \cite{cardellac14}. 
Alternatively, one may take the point of view of a center of decision (or social planner), which decides  the strategies for all players, with the goal of minimizing the global cost to the collectivity.  This leads to the concept of Pareto or cooperative equilibrium whose asymptotic formulation is reduced to the optimal control  of McKean-Vlasov dynamics for a representative player. More precisely, given the symmetry of the set-up, when the social planner  chooses  the same control policy for all the players in feedback form: $\tilde\alpha_t^i$ $=$ 
$\tilde\alpha(t,X_t^i,\bar\rho_t^N)$, $i$ $=$ $1,\ldots,N$, for some deterministic function $\tilde\alpha$ depending upon time, private state of player, 
and the empirical distribution of all players, then the theory of pro\-pagation of chaos implies that, in the limit $N$ $\rightarrow$ $\infty$, the particles $X^i$ become asymptotically independent conditionally on the random environment $W^0$, and the empirical measure $\bar\rho_t^N$ converge to the distribution $\P_{X_t}^{W^0}$ of $X_t$ given $W^0$,  and  $X$ is governed by the (stochastic) McKean-Vlasov equation: 
\beqs
dX_t &=& \tilde b(X_t,\P_{X_t}^{W^0},\tilde\alpha(t,X_t,\P_{X_t}^{W^0})) dt + \tilde \sigma(X_t,\P_{X_t}^{W^0},\tilde\alpha(t,X_t,\P_{X_t}^{W^0})) dB_t \\
& & \;\;\;\;\; + \;   \tilde \sigma_0(X_t,\P_{X_t}^{W^0},\tilde\alpha(t,X_t,\P_{X_t}^{W^0})) dW_t^0,
\enqs
for some Brownian motion $B$ independent of $W^0$.  The objective of the representative player for the Pareto equilibrium becomes the minimization of the functional
\beqs
J(\tilde\alpha) &=& \E \Big[ \int_0^T \tilde f(X_t,\P_{X_t}^{W^0},\tilde\alpha(t,X_t,\P_{X_t}^{W^0})) dt + g(X_T,\P_{X_T}^{W^0}) \Big]
\enqs
over the class of feedback controls $\tilde\alpha$. We refer to  \cite{cardel13}  
for a detailed discussion of the differences between the nature and solutions to the MFG and optimal control of McKean-Vlasov dynamics related respectively to the notions of Nash and Pareto equilibrium.  Notice that in this McKean-Vlasov control formulation,  the control $\tilde\alpha$ is of feedback (also called closed-loop) form both w.r.t. the state process $X_t$, and its conditional law process $\P_{X_t}^{W^0}$, which is $\F^0$-adapted. More generally, we can consider semi-feedback control $\alpha(t,x,\omega^0)$, in the sense that it is of closed-loop form w.r.t. the state process $X_t$, but of open-loop form w.r.t. the common noise $W^0$. In other words, one can consider  random field  control $\F^0$-progressive control process $\alpha$ $=$ 
$\{\alpha_t(x),x\in\R^d\}$, which may be viewed equivalently as processes valued in some functional space $\bA$ on $\R^d$,  typically a closed subset of the Polish space 
$C(\R^d,A)$, of continuous functions from $\R^d$ into some Euclidian space $A$. In this case, we are in the framework \reff{stoMcKean}-\reff{costMcKean} with  
$b(x,\mu,a)$ $=$ $\tilde b(x,\mu,a(x))$, $\sigma(x,\mu,a)$ $=$ $\tilde \sigma(x,\mu,a(x))$, 
$\sigma_0(x,\mu,a)$ $=$ $\tilde \sigma_0(x,\mu,a(x))$, $f(x,\mu,a)$ $=$ $\tilde f(x,\mu,a(x))$, for $(x,\mu,a)$ $\in$ 
$\R^d\times\Pc(\R^d)\times\bA$.

\vspace{1mm}
 
We also mention that  partial observation control problem arises as a particular case of our stochastic control framework \reff{stoMcKean}-\reff{costMcKean}:  Indeed, let us consider a controlled process with dynamics
\beqs
d\bar X_t &=& \bar b(\bar X_t,\alpha_t) dt + \bar \sigma(\bar X_t,\alpha_t) dB_t + \bar \sigma_0(\bar X_t,\alpha_t) dB_t^0,
\enqs
where $B,B^0$ are two independent Brownian motions on some physical probability space  $(\Omega,\Fc,\Q)$, and the signal control process can only be observed through $W^0$  given by 
\beqs
dW_t^0 &=& h(\bar X_t) dt + dB_t^0.
\enqs
The control process $\alpha$  is progressively measurable w.r.t. the observation filtration $\F^0$ gene\-rated by $W^0$,  valued typically in some Euclidian space $A$, and the cost functional to minimize over $\alpha$ is
\beqs 
J(\alpha) &=& \E^{\Q} \Big[ \int_0^T \bar f(\bar X_t,\alpha_t) dt + \bar g(\bar X_T) \Big].
\enqs
By considering the process $Z$ via
\beqs
Z_t^{-1} &=& \exp\big( - \int_0^t h(\bar X_s) dB_s^0 - \frac{1}{2} \int_0^t |h(\bar X_s)|^2 ds \big),  \; 0 \leq t \leq T, 
\enqs
the process $Z^{-1}$ is (under suitable integrability conditions on $h$) a martingale under $\Q$, and by Girsanov's theorem, this defines a probability measure $\P(d\omega)$ $=$ $Z_T^{-1}(\omega)\Q(d\omega)$, called reference probability measure, under which  the pair $(B,W^0)$ is a Brownian motion. We then see that the partial observation control problem can be recast into the framework \reff{stoMcKean}-\reff{costMcKean} of a  particular stochastic McKean-Vlasov control problem with $X$ $=$ $(\bar X,Z)$ governed by 
\beqs
d\bar X_t &=& \big( \bar b(\bar X_t,\alpha_t) - \bar\sigma_0(\bar X_t,\alpha_t) h(\bar X_t) \big) dt + \bar \sigma(\bar X_t,\alpha_t) dB_t + \bar \sigma_0(\bar X_t,\alpha_t) dW_t^0, \\
dZ_t &=& Z_t h(\bar X_t) dW_t^0,
\enqs
and a cost functional rewritten under the reference probability measure from Bayes formula as
\beqs
J(\alpha) &=& \E\Big[ \int_0^T Z_t \bar f(\bar X_t,\alpha_t) dt + Z_T \bar g(\bar X_T) \Big].
\enqs

The optimal control of McKean-Vlasov dynamics  is a rather new problem with an increasing interest  in the field of stochastic control problem. It has been studied by maximum principle methods in 
\cite{anddje10}, \cite{bucetal11}, \cite{cardel14}  for state dynamics depending upon marginal  distribution, and in \cite{carzhu14}, \cite{buclima15} for conditional McKean-Vlasov dynamics. 
This leads to a characterization of the solution in terms of an adjoint backward stochastic differential equation (BSDE) coupled with a forward SDE, and we refer to \cite{chacridel15} for a theory of  BSDE of McKean-Vlasov type.  
Alternatively, dynamic programming approach for the control of  McKean-Vlasov dynamics  has been considered in  \cite{benetal15}, \cite{benetal15b}, \cite{laupir14}  for specific McKean-Vlasov 
dynamics and under a density assumption on the probability law of the state process, and then analyzed in a general framework in \cite{phawei15} (without noise $W^0$), 
where the problem is reformulated into a deterministic control problem involving the marginal distribution process.   
%Most of these cited papers (except \cite{benetal15}) consider McKean-Vlasov equation without noise $W^0$, hence with coefficients depending on the  (deterministic) marginal  distribution.  

The aim of this paper is to develop the dynamic programming method for  stochastic McKean-Vlasov equation in a general setting. For this purpose, a key step is to show the flow property of the conditional distribution $\P_{X_t}^{W^0}$ 
of the controlled state process $X_t$ given the noise $W^0$. Then, by reformulating the original control problem into a stochastic  control problem  where  the conditional  law $\P_{X_t}^{W^0}$  is the sole controlled state variable driven by the random noise  $W^0$, and by showing the continuity of the value function in the Wasserstein space of  probability measures, we are able to prove a dynamic programming principle (DPP) for our stochastic McKean-Vlasov control problem.  Next, for exploiting the DPP, we use a notion of differentiability with respect to probability measures introduced by P.L. Lions in his lectures at the Coll\`ege de France \cite{lio12}, and detailed in the notes  \cite{car12}. 
This notion of derivative is based on the lifting of functions defined on the Hilbert space of square integrable random variables distributed according to the ``lifted" probability measure. 
By combining with a special It\^o's chain rule for flows of  conditional distributions, we derive the dynamic programming Bellman equation  for  stochastic McKean-Vlasov control problem, which is a fully nonlinear second order partial differential equation (PDE)  in the infinite dimensional Wasserstein space  of probability measures.  By adapting standard arguments to our context, we prove the viscosity property of the value function to the Bellman equation from the dynamic programming principle. To complete our PDE characterization of the value function with a uniqueness result, it is convenient to work in the lifted Hilbert space of square integrable random variables instead of the Wasserstein metric space of probability measures, in order to rely on the general results for viscosity solutions  of second order Hamilton-Jacobi-Bellman equations in separable Hilbert spaces, see \cite{lio88}, \cite{lio89b}, \cite{fabgozswi15}.  
We also state a verification theorem which is useful for getting an analytic feedback form of the optimal control when there is a smooth solution to the Bellman equation. Finally, we apply our results to the class of linear-quadratic (LQ) stochastic McKean-Vlasov control problem for which one can obtain explicit solutions, and we illustrate with an example arising from an interbank systemic risk model. 
   
The outline of the paper is organized as follows. Section 2 formulates the stochastic McKean-Vlasov control problem, and fix the standing assumptions. Section 3 is devoted to the proof and statement of the 
dynamic programming principle. We prove in Section 4 the viscosity characterization of  the value function to the Bellman equation, and the last Section 5  presents the application to the LQ framework with explicit solutions.

\section{Conditional McKean-Vlasov control problem} \label{seccondMcKean}

\setcounter{equation}{0} \setcounter{Assumption}{0}
\setcounter{Theorem}{0} \setcounter{Proposition}{0}
\setcounter{Corollary}{0} \setcounter{Lemma}{0}
\setcounter{Definition}{0} \setcounter{Remark}{0}

Let us fix some complete probability space $(\Omega,\Fc,\P)$ assumed of the form  $(\Omega^0\times\Omega^1,\Fc^0\otimes\Fc^1,\P^0\otimes\P^1)$, where $(\Omega^0,\Fc^0,\P^0)$ supports a $m$-dimensional Brownian motion $W^0$, and  $(\Omega^1,\Fc^1,\P^1)$ 
supports a $n$-dimensional Brownian motion $B$. So an element $\omega$ $\in$ $\Omega$ is written as 
$\omega$ $=$ $(\omega^0,\omega^1)$ $\in$ $\Omega^0\times\Omega^1$, and we extend canonically $W^0$ and $W$ on $\Omega$ by setting 
$W^0(\omega^0,\omega^1)$ $:=$ $W^0(\omega^0)$, $W(\omega^0,\omega^1)$ $:=$ $W(\omega^1)$, and extend similarly on $\Omega$ 
any random variable on  $\Omega^0$ or $\Omega^1$.  We assume that $(\Omega^1,\Fc^1,\P^1)$ is in the form $\Omega^1$ $=$ $\tilde\Omega^1\times\Omega^{'1}$, 
$\Fc^1$ $=$ $\Gc\otimes\Fc^{'1}$, $\P^1$ $=$ $\tilde\P^1\otimes\P^{'1}$, where $\tilde\Omega^1$ is a Polish space, $\Gc$ its Borel $\sigma$-algebra,  $\tilde\P^1$ an atomless probability  measure on 
$(\tilde\Omega^1,\Gc)$, while $(\Omega^{'1},\Fc^{'1},\P^{'1})$ supports $B$.    
We denote by  $\E^0$ (resp. $\E^1$ and $\tilde\E^1$)  the expectation under $\P^0$ (resp. $\P^1$ and $\tilde\P^1$),  by $\F^0$ $=$  
$(\Fc_t^0)_{t\geq 0}$ the $\P^0$-completion of the natural filtration generated by $W^0$ 
(and w.l.o.g. we assume that $\Fc^0$ $=$ $\Fc^0_\infty$), and by 
$\F$ $=$ $(\Fc_t)_{t\geq 0}$  the natural filtration generated by $W^0,B$, augmented with the independent $\sigma$-algebra $\Gc$.   
%in the form $\Gc$ $=$ $\{\Omega^0\times E: E \in \Gc^1\}$, for some $\sigma$-algebra $\Gc^1$ $\subset$ $\Fc^1$. 
We denote by $\Pc_{_2}(\R^d)$  the set probability measures  
$\mu$ on $\R^d$, which are square integrable, i.e. $\|\mu\|_{_2}^2$ $:=$ $\int_{\R^d} |x|^2 \mu(dx)$ $<$ $\infty$. For any $\mu$ $\in$ $\Pc_{_2}(\R^d)$,  
we denote by $L_\mu^2(\R^q)$ the set of measurable functions $\varphi$ $:$ $\R^d$ $\rightarrow$ $\R^q$, which are square integrable with respect to $\mu$, by $L_{\mu\otimes\mu}^2(\R^q)$ the set of measurable functions $\psi$ $:$ 
$\R^d\times\R^d$ $\rightarrow$ $\R^q$,  which are square integrable with respect to the product measure  $\mu\otimes\mu$, and we set
\beqs
\mu(\varphi) \; := \;   \int_{\R^d}\varphi(x)\,\mu(dx), & & \mu\otimes\mu(\psi) \;  := \;  \int_{\R^d\times\R^d} \psi(x,x') \mu(dx)\mu(dx'). 
\enqs 
We also define $L_\mu^\infty(\R^q)$ (resp.  $L_{\mu\otimes\mu}^\infty(\R^q)$) as the subset of elements $\varphi$ $\in$ $L_\mu^2(\R^q)$ (resp.  $L_{\mu\otimes\mu}^2(\R^q)$) which are bounded 
$\mu$ (resp. $\mu\otimes\mu$) a.e., and $\|\varphi\|_\infty$ is their essential supremum. 
We denote  by $L^2(\Gc;\R^d)$ (resp. $L^2(\Fc_t;\R^d)$)  the set of $\R^d$-valued  square integrable random variables on $(\tilde\Omega^1,\Gc,\tilde\P^1)$ (resp. on $(\Omega,\Fc_t,\P)$). 
For any random variable $X$ on $(\Omega,\Fc,\P)$, we denote by $\P_X$  its probability law (or distribution) under  $\P$,  and we know that 
$\Pc_{_2}(\R^d)$ $=$ $\{\P_\xi = \tilde\P^1_\xi: \xi \in L^2(\Gc;\R^d)\}$ since $(\tilde\Omega^1,\Gc,\tilde\P^1)$ is Polish and atomless (we say that $\Gc$ is rich enough). 
We  often   write $\Lc(\xi)$ $=$ $\P_\xi$ $=$ $\tilde\P^1_\xi$ for the law of 
$\xi$ $\in$ $L^2(\Gc;\R^d)$.  
The space  $\Pc_{_2}(\R^d)$ is a metric space equipped with the $2$-Wasserstein distance
\beqs
\Wc_2(\mu,\mu') &:=& \inf\Big\{ \Big( \int_{\R^d\times\R^d} |x-y|^2 \pi(dx,dy)\Big)^{1\over 2}:
\pi \in \Pc_{_2}(\R^d\times\R^d) \mbox{ with marginals } \mu \mbox{ and } \mu' \Big\} \\
&=& \inf\Big\{  \Big(\E|\xi-\xi'|^2\Big)^{1\over 2}: \;\; \xi,\xi' \in L^2(\Gc;\R^d) \mbox{ with } \Lc(\xi) = \mu, \; \Lc(\xi') = \mu' \Big\}, 
\enqs
and endowed with the corresponding Borel $\sigma$-field $\Bc(\Pc_{_2}(\R^d))$.  We recall in the next remark some useful properties on this Borel 
$\sigma$-field.

\begin{Remark} \label{remBorel}
{\rm Denote by $\mathscr C_{_2}(\R^d)$ the set of continuous  functions on $\R^d$ with quadratic growth, and  for any  
$\varphi\in\mathscr C_{_2}(\R^d)$, define the map $\Lambda_{_\varphi}\colon\Pc_{_2}(\R^d)\rightarrow\R$ by $\Lambda_{_\varphi}\mu$ $=$ 
$\mu(\varphi)$, for $\mu$ $\in$ $\Pc_{_2}(\R^d)$. 
%as follows:
%\beqs
%\Lambda_{_\varphi} \mu  & := & \mu(\phi) \; = \;  \int_{\R^d}\varphi(x)\,\mu(dx), \qquad \text{for every }\pi\in\Pc_{_2}(\R^d).
%\enqs
By Theorem 7.12 in \cite{vil03}, for $(\mu_n)_n$, $\mu$ $\in$ $\Pc_{_2}(\R^d)$, we have that $\Wc_2(\mu_n,\mu)\rightarrow0$ if and only if, for every 
$\varphi\in\mathscr C_{_2}(\R^d)$, $\Lambda_{_\varphi}\mu_n\rightarrow\Lambda_{_\varphi}\mu$. 
Therefore, recalling also that $(\Pc_{_2}(\R^d),\Wc_2)$ is a complete separable metric space (see e.g. Proposition 7.1.5 in \cite{ambetal05}), we notice that  $\Bc(\Pc_{_2}(\R^d))$ coincides with the cylindrical  $\sigma$-algebra $\sigma(\Lambda_{_\varphi},\,\varphi\in\mathscr C_{_2}(\R^d))$.   
Consequently,  given a mea\-surable space $(E,\Ec)$ and a map $\rho$ $:$ $E$ $\rightarrow$ $\Pc_{_2}(\R^d)$, $\rho$ is measurable if and only if  the map  $\Lambda_{_\varphi} \circ \rho$ $=$ $\rho(\varphi)$ $:$ $E$ $\rightarrow$ $\R$  is measurable, for any $\varphi\in\mathscr C_{_2}(\R^d)$. 
%which means that the map $\Lambda_\varphi$ is $\Bc(\Pc_{_2}(\R^d))$-measurable, 
%for any $\varphi\in\mathscr C_{_2}(\R^d)$, and also for any measu\-rable function $\varphi$ with quadratic growth condition, by a monotone class argument. %Consequently,  given a measurable space $(E,\Ec)$ and a map $\rho$ $:$ $E$ $\rightarrow$ $\Pc_{_2}(\R^d)$, $\rho$ is measurable if and only if  the map  
%$\Lambda_{_\varphi} \circ \rho$ $=$ $\rho(\varphi)$ $:$ $E$ $\rightarrow$ $\R$  is measurable, for any $\varphi\in\mathscr C_{_2}(\R^d)$.  
Finally, we notice that the map $\Lambda_\varphi$ is $\Bc(\Pc_{_2}(\R^d))$-measurable, for any measurable function 
$\varphi$ with quadratic growth condition, by using a monotone class argument since it holds true whenever $\varphi\in\mathscr C_{_2}(\R^d)$. 
}
\ep
\end{Remark}

\vspace{1mm}

%coincides with the $\sigma$-field generated by the mappings $\mu$ $\in$ $\Pc_{_2}(\R^d)$ $\mapsto$ $\mu(B)$, $B$ running  over the Borel subsets of $\R^d$, see e.g.  chapter 7 in \cite{vil03}.   

\noindent  $\bullet$ {\it Admissible controls}.  We are given a Polish set $\bA$ equipped with the distance $d_{_A}$, satisfying w.l.o.g. $d_{_A}$ $<$ $1$,  representing the control set, and  
we denote by $\Ac$ the set of $\F^0$-progressive processes $\alpha$ valued in $\bA$. 
%and s.t. $\E^0\big[ \int_0^T \|\alpha_t\|_{_A}^2 dt \big]$  $<$ $\infty$. 
Notice that $\Ac$ is a separable metric  space endowed with the Krylov distance 
$\Delta(\alpha,\beta)$ $=$ $\E^0[\int_0^T d_{_A}(\alpha_t,\beta_t) dt]$.  We denote by $\Bc_{\Ac}$ the Borel $\sigma$-algebra of $\Ac$.

\vspace{2mm} 

\noindent $\bullet$ {\it  Controlled stochastic McKean-Vlasov dynamics}.  For $(t,\xi)$ $\in$ $[0,T]\times L^2(\Fc_t;\R^d)$,  and given 
 $\alpha$ $\in$ $\Ac$,  we consider the stochastic  McKean-Vlasov equation: 
\begin{equation} \label{Mckean}
\left\{
\begin{array}{rcl}
dX_s &=&  b(X_s,\P_{X_s}^{W^0},\alpha_s) ds + \sigma(X_s,\P_{X_s}^{W^0},\alpha_s) dB_s \\
& & \;\;\;\;\;\;\;\;\;\; + \;  \sigma_0(X_s,\P_{X_s}^{W^0},\alpha_s) dW_s^0, \;\;\;  t \leq s \leq T, \\
X_t &=& \xi. 
\end{array}
\right.
\end{equation}
Here, $\P_{X_s}^{W^0}$ denotes the regular conditional distribution of $X_s$ given $\Fc^0$, and its realization at some $\omega^0$ $\in$ $\Omega^0$ also reads as the law under $\P^1$ of the random variable $X_s(\omega^0,.)$ on 
$(\Omega^1,\Fc^1,\P^1)$, i.e. $\P_{X_s}^{W^0}(\omega^0)$ $=$ $\P^1_{X_s(\omega^0,.)}$. 
The coefficients $b$, $\sigma$, $\sigma_0$ are measurable functions from $\R^d\times\Pc_{_2}(\R^d)\times \bA$ into $\R^d$, respectively 
$\R^{d\times n}$, $\R^{d\times m}$, and satisfy the condition: 

\vspace{2mm}

\noindent {\bf (H1)} 
\begin{itemize}
\item[(i)] There exists some positive constant $C$ s.t. for all  $x,x'$ $\in$ $\R^d$, $\mu,\mu'$ $\in$ $\times\Pc_{_2}(\R^d)$, and $a$ $\in$  $\bA$, 
\beqs
& & |b(x,\mu,a) - b(x',\mu',a)| + |\sigma(x,\mu,a) - \sigma(x',\mu',a)| + |\sigma_0(x,\mu,a) - \sigma_0(x',\mu',a)|  \\
& \leq & C \Big( |x-x'|  + \Wc_2(\mu,\mu') \Big), 
\enqs
and 
\beqs
|b(0,\delta_0,a)| +  |\sigma(0,\delta_0,a)| + |\sigma_0(0,\delta_0,a)| & \leq & C.  
\enqs 
\item[(ii)] For all $(x,\mu)$ $\in$ $\R^d\times\Pc_{_2}(\R^d)$, the functions $a$  $\mapsto$ $b(x,\mu,a)$, $\sigma(x,\mu,a)$, $\sigma_0(x,\mu,a)$ are continuous on $\bA$.  
\end{itemize}

\begin{Remark}
{\rm  We have chosen a control formulation where the process $\alpha$ is required to be progressively measurable w.r.t. the filtration $\F^0$ of the sole common noise. 
%hence viewed as a semi-closed loop control when ${\bf A}$ is a functional space. 
This form is used for rewriting  the cost functional in terms of the conditional law as sole state variable, see \reff{Jmu}, which is then convenient for deriving the dynamic programming principle. 
In the case where ${\bf A}$ is a functional space on the state space $\R^d$, meaning that $\alpha$ is a semi closed-loop control,  and when the coefficients are in the form: 
$b(x,\mu,a)$ $=$ $\tilde b(x,\mu,a(x))$, $\sigma(x,\mu,a)$ $=$ $\tilde \sigma(x,\mu,a(x))$, $\sigma_0(x,\mu,a)$ $=$ $\tilde \sigma_0(x,\mu,a(x))$ (see discussion in the introduction), the Lipschitz condition in 
{\bf (H1)}(i) requires that $a$ $\in$ ${\bf A}$ is Lipschitz continuous with a prescribed Lipschitz constant, which is somewhat a restrictive condition.  The more general case where the control $\alpha$ is allowed to be measurable 
with respect to the filtration $\F$ of  both noises, i.e.,  $\alpha$ of open-loop form, is certainly an important extension, and  left for future work. In this case, one should consider as state variables the pair composed of the process 
$X_t$ and its conditional law $\P_{_{X_t}}^{W^0}$, see the recent paper  \cite{BCP16} where a dynamic programming principle is stated when the control is allowed to be of open-loop form in the case without common noise. 
}
\ep
\end{Remark}

Under {\bf (H1)}(i), there exists a unique solution to \reff{Mckean} (see e.g. \cite{kurxio99}), 
denoted by $\{X_s^{t,\xi,\alpha},t\leq s \leq T\}$, which is $\F$-adapted, and satisfies the square-integrability condition: 
\beq \label{momentX}
\E \Big[ \sup_{t\leq s \leq T} |X_s^{t,\xi,\alpha}|^2 \Big] & \leq & C \Big( 1  + \E|\xi|^2  \Big) \; < \; \infty, 
\enq
for some positive constant $C$ independent of $\alpha$. 
We shall sometimes omit the dependence of  $X^{t,\xi}$ $=$ $X^{t,\xi,\alpha}$ on $\alpha$ when there is no ambiguity. 
Since $\{X_s^{t,\xi},t\leq s \leq T\}$ is $\F$-adapted, and  $W^0$ is a $(\P,\F)$-Wiener process, we notice that 
$\P_{X_s^{t,\xi}}^{W^0}(dx)$ $=$ $\P[X_s^{t,\xi} \in dx | \Fc^0]$ $=$ $\P[X_s^{t,\xi} \in dx | \Fc_s^0]$. We thus have for any $\varphi$ $\in$ $\mathscr C_{_2}(\R^d)$:  
\beq \label{defPW0}
\P_{X_s^{t,\xi}}^{W^0}(\varphi) &=& \E \Big[ \varphi(X_s^{t,\xi}) \big| \Fc^0 \Big] \; = \; \E \Big[ \varphi(X_s^{t,\xi}) \big| \Fc_s^0 \Big], \;\;\; t \leq s \leq T,
\enq
which shows that $\P_{X_s^{t,\xi}}^{W^0}(\varphi)$ is $\Fc_s^0$-measurable, and therefore, in view of the measurability property in Remark \ref{remBorel}, that  
$\{\P_{X_s^{t,\xi}}^{W^0}, t\leq s \leq T\}$ is $(\Fc_s^0)_{t\leq s\leq T}$-adapted.  Moreover, since $\{\P_{X_s^{t,\xi}}^{W^0}, t\leq s \leq T\}$ is valued in $\Pc_{_2}(C([t,T];\R^d))$, the set of square integrable probability measures on the space $C([t,T];\R^d)$ of continuous functions from $[t,T]$ into $\R^d$, it also has continuous trajectories, and is then $\F^0$-progressively measurable (actually even $\F^0$-predictable).

\vspace{1mm} 

\noindent $\bullet$  {\it Cost functional and value function}. We are given a running cost function $f$ defined on $\R^d\times\Pc_{_2}(\R^d)\times\bA$, and 
a terminal cost function $g$ defined on $\R^d\times\Pc_{_2}(\R^d)$,  assumed to satisfy the condition

\vspace{1mm}

\noindent {\bf (H2)} 
\begin{itemize}
\item[(i)]
There exists some positive constant $C$ s.t. for all $(x,\mu,a)$ $\in$ $\R^d\times\Pc_{_2}(\R^d)\times\bA$,
\beqs
|f(x,\mu,a)| + |g(x,\mu)| & \leq & C \big(1 + |x|^2  + \|\mu\|_{_2}^2 \big). 
\enqs
\item[(ii)] The functions $f$, $g$ are continuous on $\R^d\times\Pc_{_2}(\R^d)\times\bA$, resp. on  $\R^d\times\Pc_{_2}(\R^d)$,  and satisfy the local Lipschitz condition, uniformly w.r.t. $\bA$: there exists some positive constant 
$C$ s.t. for all $x,x'$ $\in$ $\R^d$, $\mu,\mu'$ $\in$ $\Pc_{_2}(\R^d)$,  $a$ $\in$ $\bA$, 
\beqs
& & | f(x,\mu,a) - f(x',\mu',a)| + |g(x,\mu) - g(x',\mu')| \\
& \leq & C(1 +|x| + |x'| + \|\mu\|_{_2} + \|\mu'\|_{_2} )  \big(|x-x'| + \Wc_2(\mu,\mu') \big). 
\enqs
\end{itemize}

We then consider the cost functional: 
\beqs
J(t,\xi,\alpha) & :=& \E \Big[ \int_t^T f\big(X_s^{t,\xi},\P_{X_s^{t,\xi}}^{W^0},\alpha_s \big) ds 
+ g\big(X_T^{t,\xi},\P_{X_T^{t,\xi}}^{W^0}\big) \Big], 
\enqs
which is well-defined and finite for all $(t,\xi,\alpha)$ $\in$ $[0,T]\times L^2(\Gc;\R^d)\times\Ac$, and we define the value function of the 
conditional McKean-Vlasov control problem as
\beq \label{defv}
v(t,\xi) &:=& \inf_{\alpha\in\Ac} J(t,\xi,\alpha), \;\;\;\;\;   (t,\xi) \in [0,T]\times L^2(\Gc;\R^d). 
\enq
From the estimate \reff{momentX} and the growth condition in {\bf (H2)}(i), it is clear that $v$ also satisfies a quadratic growth condition:
\beq \label{vquadra}
|v(t,\xi)| & \leq & C\big(1 + \E|\xi|^2 \big), \;\;\; \forall \xi \in  L^2(\Gc;\R^d). 
\enq

Our goal is to characterize the value function $v$ as solution of a partial differential equation by means of a dynamic programming approach. 
%In the rest of the paper, we make the standing assumptions {\bf (H1)}-{\bf (H2)}. 

\section{Dynamic programming}

\setcounter{equation}{0} \setcounter{Assumption}{0}
\setcounter{Theorem}{0} \setcounter{Proposition}{0}
\setcounter{Corollary}{0} \setcounter{Lemma}{0}
\setcounter{Definition}{0} \setcounter{Remark}{0}

%\subsection{Dynamic programming principle}

The aim  of this section is to prove  the dynamic programming principle (DPP)  for the value function $v$ in \reff{defv}  of the conditional McKean-Vlasov control problem.

\subsection{Flow properties}

We shall assume that $(\Omega^0,W^0,\P^0)$ is the canonical space, i.e. $\Omega^0$ $=$ $C(\R_+,\R^m)$, the set of continuous functions from $\R_+$ into $\R^m$, 
$W^0$ is the canonical process,  and $\P^0$ the Wiener measure. 
Following \cite{clataltan15}, we introduce the class of shifted control processes constructed by concatenation of paths: for $\alpha$ $\in$ $\Ac$, $(t,\bar\omega^0)$ $\in$ $[0,T]\times\Omega^0$, we set
\beqs
\alpha_s^{t,\bar\omega^0}(\omega^0) &:=& \alpha_s(\bar\omega^0\otimes_t\omega^0),  \;\;\; (s,\omega^0) \in [0,T]\times\Omega^0, 
\enqs
where $\bar\omega^0\otimes_t\omega^0$ is the element in $\Omega^0$ defined by
\beqs
\bar\omega^0\otimes_t\omega^0(s) &:=& \bar\omega^0(s) 1_{s < t} + \big( \bar\omega^0(t) +  \omega^0(s)-\omega^0(t) \big) 1_{s\geq t}. 
\enqs
We notice that for fixed $(t,\bar\omega^0)$, the process $\alpha^{t,\bar\omega^0}$ lies in $\Ac_t$, the set of  elements in $\Ac$ which are independent of $\Fc_t^0$ under $\P^0$.  
For any $\alpha$ $\in$  $\Ac$, and $\F^0$-stopping time $\theta$, we denote by $\alpha^\theta$ the  map 
\beqs
\alpha^\theta: (\Omega^0,\Fc_\theta^0) & \rightarrow & (\Ac,\Bc_\Ac) \\
\omega^0 & \mapsto & \alpha^{\theta(\omega^0),\omega^0}. 
\enqs
%process defined on $(\Omega^0,\Fc^0,\P^0)$ by 
%\beqs
%\alpha^\theta(\omega^0) &:= & \alpha^{\theta(\omega^0),\omega^0}(\omega^0). 
%\enqs

The key step in the proof the DPP is to obtain a flow property on the controlled conditional distribution  $\F^0$-progressively measurable process $\{\P_{X_s^{t,\xi}}^{W^0}, t\leq s \leq T\}$, for  
$(t,\xi)$ $\in$ $[0,T]\times L^2(\Fc_t;\R^d)$, and $\alpha$ $\in$ $\Ac$.  
%We denote by $L^2(\Fc_t^0;\Pc_{_2}(\R^d))$ the set of $\Fc_t^0$-measurable probability measures $\mu_t$ valued in $\Pc_{_2}(\R^d)$, which are square integrable, i.e.  $\E^0[ \|\mu_t\|_{_2}^2]$ $<$ $\infty$.
%and by observing that $\E^0[ \|\P_{X_s^{t,\xi}}^{W^0}\|_{_2}^2]$ $=$ $\E[|X_s^{t,\xi}|^2]$, we notice that  $\P_{X_s^{t,\xi}}^{W^0}$ lies in 
%$L^2(\Fc_s^0;\Pc_{_2}(\R^d))$ for  $t\leq s\leq T$. 

\begin{Lemma} \label{lemrho}
For any $t$ $\in$ $[0,T]$, $\mu$ $\in$ $\Pc_{_2}(\R^d)$, $\alpha$ $\in$ $\Ac$,  the relation  given by 
\beq \label{defrho}
\rho_s^{t,\mu,\alpha} &:=&  \P_{X_s^{t,\xi,\alpha}}^{W^0}, \;\;\; t \leq s \leq T, \; \mbox{ for } \; \xi \in L^2(\Fc_t;\R^d) \; \mbox{ s.t. }  \;  \P_\xi^{W^0} = \mu,
\enq
defines a square integrable  $\F^0$-progressive continuous process in $\Pc_{_2}(\R^d)$. Moreover,  the map $(s,t,\omega^0,\mu,\alpha)$ $\in$ $[0,T]\times [0,T]\times\Omega^0\times\Pc_{_2}(\R^d)\times\Ac$ $\rightarrow$ $\rho_s^{t,\mu,\alpha}(\omega^0)$ $\in$ $\Pc_{_2}(\R^d)$ (with the convention that $\rho_s^{t,\mu,\alpha}$ $=$ $\mu$ for $s$ $\leq$ $t$) is measurable,  and satisfies the flow pro\-perty: 
$\rho_s^{t,\mu,\alpha}$ $=$ $\rho_s^{\theta,\rho_\theta^{t,\mu,\alpha},\alpha^\theta}$, $\P^0$-a.s.,  i.e. 
\beq \label{flowrho}
\rho_s^{t,\mu,\alpha}(\omega^0)  &=& \rho_s^{\theta(\omega^0),\rho_{\theta(\omega^0)}^{t,\mu,\alpha}(\omega^0),\alpha^{\theta(\omega^0),\omega^0}}(\omega^0), 
%\;\;\; \mbox{ for }   (\tau,\pi,\beta) \; = \;  (\theta(\omega^0),\rho_{\tau}^{t,\mu_t,\alpha}(\omega^0),\alpha^{\theta(\omega^0),\omega^0})  
\;\;\;  s \in [\theta,T], \; \P^0(d\omega^0)-a.s
\enq
for all $\theta$ $\in$ $\Tc^0_{t,T}$, the set of  $\F^0$-stopping times valued in $[t,T]$.  
\end{Lemma}
{\bf Proof.} {\bf 1.} First observe that for any $t$ $\in$ $[0,T]$, $\xi$ $\in$ $L^2(\Fc_t;\R^d)$, $\alpha$ $\in$ $\Ac$, we have: $\E^0[ \|\P_{X_s^{t,\xi,\alpha}}^{W^0}\|_{_2}^2]$ $=$ $\E[|X_s^{t,\xi,\alpha}|^2]$ $<$ $\infty$, 
%i.e. $\P_{X_s^{t,\xi,\alpha}}^{W^0}$ $\in$  $L^2(\Fc_s^0;\Pc_{_2}(\R^d))$, 
which means that the process $\{\P_{X_s^{t,\xi,\alpha}}^{W^0},t\leq s\leq T\}$ is square integrable, and we recall (see the discussion after \reff{defPW0}) that it is $\F^0$-progressively measurable.  
\begin{itemize}
\item[(i)] Notice  that for $\P^0$-a.s $\omega^0$ $\in$ $\Omega^0$, the law of the solution $\{X_s^{t,\xi,\alpha}(\omega^0,.),t\leq s\leq T\}$ to \reff{Mckean} on 
$(\Omega^1,\Fc^1,\P^1)$ is unique in law, which implies  that $ \P_{X_s^{t,\xi,\alpha}}^{W^0}(\omega^0)$ $=$ $\P^1_{X_s^{t,\xi,\alpha}(\omega^0,.)}$, 
$t\leq s \leq T$, depends on $\xi$ only through $\P_\xi^{W^0}(\omega^0)$ $=$ $\P^1_{\xi(\omega^0,.)}$.  In other words, for any $\xi_1$, $\xi_2$ $\in$ $L^2(\Fc_t;\R^d)$ s.t. $\P_{\xi_1}^{W^0}$ $=$ $\P_{\xi_2}^{W^0}$, 
the processes $\{\P_{X_s^{t,\xi_1,\alpha}}^{W^0},t\leq s\leq T\}$ and  $\{\P_{X_s^{t,\xi_2,\alpha}}^{W^0},t\leq s\leq T\}$ are indistinguishable. 
%we have 
%\beqs
%\P_{X_s^{t,\xi_1}}^{W^0} &=&  \P_{X_s^{t,\xi_2}}^{W^0}, \;\;\; t \leq s \leq T, \;\; \P^0-a.s.
%\enqs
\item[(ii)] Let us now check that for any $\mu$ $\in$ $\Pc_{_2}(\R^d)$,  one can find $\xi$ $\in$ $L^2(\Fc_t;\R^d)$ s.t. $\P_{\xi}^{W^0}$ $=$ $\mu$.  
Indeed, recalling that  $\Gc$ is rich enough, one can find $\xi$ $\in$ $L^2(\Gc;\R^d)$ $\subset$ $L^2(\Fc_t;\R^d)$  s.t. $\Lc(\xi)$ $=$ $\mu$. Since $\Gc$ is independent of $W^0$, this also means that 
$\P_\xi^{W^0}$ $=$ $\mu$. 
%for $\P^0$-a.s $\omega^0$ $\in$ $\Omega^0$,  a random variable $\xi(\omega^0,.)$ $\in$ 
%$L^2(\Gc^1;\R^d)$ s.t. $\P^1_{\xi(\omega^0,.)}$ $=$ $\P^{W^0}_\xi(\omega^0)$ $=$ $\mu_t(\omega^0)$. Since $\mu_t$ is $\Fc_t^0$-measurable, the map $\omega^0$ $\rightarrow$ $\xi(\omega^0,.)$ can be also chosen 
%$\Fc_t^0$-measurable (BY A MEASURABLE SELECTION ARGUMENT?),  and so the random variable $\xi$ on $(\Omega,\Fc,\P)$ is $\Fc_t$-measurable. It is also square integrable, since 
%$\E|\xi|^2$ $=$ $\E^0[\|\mu_t\|_{_2}^2]$ $<$ $\infty$, and then $\xi$ $\in$ $L^2(\Fc_t;\R^d)$.  
\end{itemize}
In view of the uniqueness result in (i), and the representation result in (ii), one can define the process $\{\rho_s^{t,\mu,\alpha},t\leq s\leq T\}$ by the relation \reff{defrho}, and this process is a square integrable  $\F^0$-progressively measurable process in 
$\Pc_{_2}(\R^d)$.  
%It has also continuous trajectories as $\rho^{t,\mu,\alpha}$ is valued in $\Pc_{_2}(C([t,T];\R^d))$, the set of square integrable probability measures on the space $C([t,T];\R^d)$ of continuous functions from $[t,T]$ into $\R^d$. 

\vspace{1mm}

\noindent {\bf 2.} Let us now prove the joint measurability of  $\rho_s^{t,\mu,\alpha}(\omega^0)$ in $(t,s,\omega^0,\mu,\alpha)$ $\in$ $[0,T]\times [0,T]\times\Omega^0\times\Pc_{_2}(\R^d)\times\Ac$. 
Given $t$ $\in$ $[0,T]$, $\mu$ $\in$ $\Pc_{_2}(\R^d)$, $\alpha$ $\in$ $\Ac$, let $\xi$ $\in$ $L^2(\Gc;\R^d)$ s.t. $\Lc(\xi)$ $=$ $\mu$. We construct $X^{t,\xi,\alpha}$ using Picard's iteration by defining recursively a sequence of processes 
$(X^{(m),t,\xi,\alpha})_m$  as follows: we start from $X^{(0),t,\xi,\alpha}$ $\equiv$ $0$, and define $\rho^{(0),t,\mu,\alpha}$ by formula \reff{defrho} with $X^{(0),t,\xi,\alpha}$ instead of $X^{t,\xi,\alpha}$, and see that $\rho^{(0),t,\mu,\alpha}$ $=$ $\delta_0$. 
\begin{itemize}
\item[-]   The process  $X^{(1),t,\xi,\alpha}$ is given by
\beqs
X_s^{(1),t,\xi,\alpha} &=& \xi + \int_t^s b(0,\delta_0,\alpha_r) dr +  \int_t^s \sigma(0,\delta_0,\alpha_r) dB_r +  \int_t^s \sigma_0(0,\delta_0,\alpha_r) dW_r^0,
\enqs
for $0\leq t\leq s\leq T$  (and $X_s^{(1),t,\xi,\alpha}$ $=$ $\xi$ when $s$ $<$ $t$), and we notice that the map $X^{(1),t,\xi,\alpha}$ $:$ 
$([t,T]\times\Omega,\Bc([t,T])\otimes\Fc)$ $\rightarrow$  $(\R^d,\Bc(\R^d))$ is measurable, up to indistinguishability.   
We then define $\rho^{(1),t,\mu,\alpha}$ by formula \reff{defrho} with $X^{(1),t,\xi,\alpha}$ instead of $X^{t,\xi,\alpha}$, so that 
\beqs
\rho_s^{(1),t,\mu,\alpha}(\omega^0)(\varphi) &=& \E^1 \Big[ \varphi\big(X_s^{(1),t,\xi,\alpha}(\omega^0,.) \big) \Big]  \; = \; \int_{\R^d} \Phi^{(1)}(x,t,s,\omega^0,\alpha) \mu(dx),
\enqs
for any $\varphi\in\mathscr C_{_2}(\R^n)$, where $\Phi^{(1)}$ $:$ $\R^d\times[0,T]\times [0,T]\times\Omega^0\times\Ac$ $\rightarrow$ $\R$ is measurable with quadratic growth condition in $x$, uniformly in $(t,s,\omega^0,\alpha)$, and given by:
\beqs
\Phi^{(1)}(x,t,s,\omega^0,\alpha) &=& \E^1 \Big[ \varphi\big(x + \int_t^s b(0,\delta_0,\alpha_r(\omega^0)) dr +  \int_t^s \sigma(0,\delta_0,\alpha_r(\omega^0)) dB_r   \\
& & \hspace{2cm}  + \;   \int_t^s \sigma_0(0,\delta_0,\alpha_r(\omega^0)) dW_r^0(\omega^0) \big) \Big], \; t \leq s \leq T, 
\enqs
and $\Phi^{(1)}(x,t,s,\omega^0,\alpha)$ $=$ $\varphi(x)$ when $s$ $<$ $t$. 
By a monotone class argument (first considering the case when  $\Phi^{(1)}(x,t,s,\omega^0,\alpha)$ is expressed as a product $h(x)\ell(t,s,\omega^0,\alpha)$ for some measurable and bounded functions $h,\ell$), we deduce that 
$\rho_s^{(1),t,\mu,\alpha}(\omega^0)(\varphi)$  is jointly measurable in $(t,s,\omega^0,\mu,\alpha)$. By  Remark \ref{remBorel}, this means that  the map $(t,s,\omega^0,\mu,\alpha)$ $\in$ $[0,T]\times [0,T]\times\Omega^0\times\Pc_{_2}(\R^d)\times\Ac$ $\mapsto$ $\rho_s^{(1),t,\mu,\alpha}(\omega^0)$ $\in$ 
$\Pc_{_2}(\R^d)$ is measurable. 
\item[-] We define recursively $X^{(m+1),t,\xi,\alpha}$ assuming that $X^{(m),t,\xi,\alpha}$ has been already defined. We assume that the map 
$X^{(m),t,\xi,\alpha}$ $:$ $([t,T]\times\Omega,\Bc([t,T])\otimes\Fc)$ $\rightarrow$  $(\R^d,\Bc(\R^d))$ is measurable (up to indistinguishability), 
and we define  $\rho_s^{(m),t,\mu,\alpha}(\omega^0)$ given by formula \reff{defrho} with $X^{(m),t,\xi,\alpha}$ instead of $X^{t,\xi,\alpha}$.  Moreover, we  suppose that 
$\rho_s^{(m),t,\mu,\alpha}(\omega^0)$ is jointly measurable in $(t,s,\omega^0,\mu,\alpha)$.  Then, we define the process $X^{(m+1),t,\xi,\alpha}$ as follows:
\beqs
X_s^{(m+1),t,\xi,\alpha} &=& \xi + \int_t^s b(X_r^{(m),t,\xi,\alpha},\rho_r^{(m),t,\mu,\alpha},\alpha_r) dr +  \int_t^s \sigma(X_r^{(m),t,\xi,\alpha},\rho_r^{(m),t,\mu,\alpha},\alpha_r) dB_r \\
& & \;\;\;\;\;\;\; + \;   \int_t^s \sigma_0(X_r^{(m),t,\xi,\alpha},\rho_r^{(m),t,\mu,\alpha},\alpha_r) dW_r^0,
\enqs
for $0\leq t\leq s\leq T$ (and $X_s^{(m+1),t,\xi,\alpha}$ $=$ $\xi$ when $s$ $<$ $t$), and notice by construction that  the map
$X^{(m+1),t,\xi,\alpha}$ $:$ 
$[t,T]\times\Omega,\Bc([t,T])\otimes\Fc)$ $\rightarrow$  $(\R^d,\Bc(\R^d))$ is measurable, up to indistinguishability.
We can then define $\rho^{(m+1),t,\mu,\alpha}$ by formula \reff{defrho} with $X^{(m+1),t,\xi,\alpha}$ instead of $X^{t,\xi,\alpha}$, namely
%so that for any  $\varphi\in\mathscr C_{_2}(\R^n)$, 
\beqs
\rho_s^{(m+1),t,\mu,\alpha}(\omega^0)(\varphi) &=& \E^1 \Big[ \varphi\big(X_s^{(m+1),t,\xi,\alpha}(\omega^0,.) \big) \Big],  
%\; = \; \int_{\R^d} \Phi^{(m+1)}(x,t,s,\omega^0,\mu,\alpha) \mu(dx),
\enqs
for any $\varphi\in\mathscr C_{_2}(\R^n)$, $\omega^0$ $\in$ $\Omega^0$. From the  (iterated) dependence of $X^{(m+1),t,\xi,\alpha}$ on $\xi$, 
and by Fubini's theorem (recalling the product structure of the probability space $\Omega^1$ on which are defined the random variable $\xi$ 
of law $\mu$ and the Brownian motion $B$), we then have
\beqs
\rho_s^{(m+1),t,\mu,\alpha}(\omega^0)(\varphi) 
%&=& \E^1\Big[ \varphi\Big( \xi + \int_t^s b(\xi + \ldots,\rho_r^{(m),t,\mu,\alpha},\alpha_r) dr \\
%& & \;\;\;\;\;\;\;\;\;\; + \;  \int_t^s \sigma(\xi + \ldots,\rho_r^{(m),t,\mu,\alpha},\alpha_r) dB_r \\
%& & \;\;\;\;\;\;\;\;\;\;  + \;  \int_t^s \sigma_0(\xi + \ldots,\rho_r^{(m),t,\mu,\alpha},\alpha_r) dW_r(\omega^0) \Big) \Big] \\
&=& \int_{\R^d} \Phi^{(m+1)}(x,t,s,\omega^0,\mu,\alpha) \mu(dx),
\enqs
%by Fubini's theorem (recalling the product structure of the probability space $\Omega^1$ on which are defined the random variable $\xi$ of law $\mu$ 
%and the Brownian motion $B$), and 
where $\Phi^{(m+1)}$ $:$ $\R^d\times[0,T]\times [0,T]\times\Omega^0\times\Pc_{_2}(\R^d)\times\Ac$ 
$\rightarrow$ $\R$ is measurable with quadratic growth condition uniformly in $(t,s,\omega^0,\alpha)$, and  given by
\beqs
\Phi^{(m+1)}(x,t,s,\omega^0,\mu,\alpha) &=& \E^1 \Big[ \varphi\Big( x + \int_t^s b(x + \ldots,\rho_r^{(m),t,\mu,\alpha},\alpha_r) dr \\
& & \;\;\;\;\;\;\;\; + \;  \int_t^s \sigma(x + \ldots,\rho_r^{(m),t,\mu,\alpha},\alpha_r) dB_r \\
& & \;\;\;\;\;\;\;\;  + \;  \int_t^s \sigma_0(x + \ldots,\rho_r^{(m),t,\mu,\alpha},\alpha_r) dW_r(\omega^0) \Big) \Big], \; t \leq s \leq T, 
\enqs
and $\Phi^{(m+1)}(x,t,s,\omega^0,\alpha)$ $=$ $\varphi(x)$ when $s$ $<$ $t$.   
We then see that $\rho_s^{(m+1),t,\mu,\alpha}(\omega^0)(\varphi)$ is jointly measurable in 
$(t,s,\omega^0,\mu,\alpha)$  (using again a monotone class argument), and deduce by Remark \ref{remBorel} that  the map $(t,s,\omega^0,\mu,\alpha)$ $\in$ $[0,T]\times [0,T]\times\Omega^0\times\Pc_{_2}(\R^d)\times\Ac$ $\mapsto$ 
$\rho_s^{(m+1),t,\mu,\alpha}(\omega^0)$ $\in$  $\Pc_{_2}(\R^d)$ is measurable.  
\end{itemize} 
Now that we have constructed the sequence $(X^{(m),t,\xi,\alpha})_m$, one can show by proceeding along the same lines as in the proof of Theorem IX.2.1 in \cite{revyor99} or Theorem V.8 in \cite{pro04} that
\beqs
\sup_{t\leq s\leq T} |X_s^{(m),t,\xi,\alpha} - X_s^{t,\xi,\alpha} | & \underset{m\rightarrow\infty}{\overset{\P}{\longrightarrow}}  & 0, 
\enqs
where the convergence holds in probability. Then, by the same arguments as in the proof of Lemma 3.2 in \cite{BCP16} (see their Appendix B), this implies that  the following convergence holds in probability:  
\beqs
\Wc_{\text{\tiny$2$}}\big(\rho_s^{(m),t,\mu,\alpha},\rho_s^{t,\mu,\alpha}\big) & \underset{m \rightarrow\infty}{\overset{\P^0}{\longrightarrow}} & 0, 
\enqs
for all $s$ $\in$ $[t,T]$, $\mu$ $\in$ $\Pc_{_2}(\R^d)$, and $\alpha$ $\in$ $\Ac$. Since for any $m$ $\in$ $\N$, $\rho_s^{(m),t,\mu,\alpha}(\omega^0)$ is jointly measurable in $(t,s,\omega^0,\mu,\alpha)$, we deduce by  
proceeding for instance as in the first item of Exercise IV.5.17 in \cite{revyor99}, 
and recalling that $\Fc^0$ is assumed to be a complete $\sigma$-field, that 
%there is a \rd{version} of  $(\rho_s^{t,\mu,\alpha})_{t\leq s\leq T}$ (still denoted by $(\rho_s^{t,\mu,\alpha})_{t\leq s\leq T})$, such that 
the map $(t,s,\omega^0,\mu,\alpha)$ $\in$ $[0,T]\times [0,T]\times\Omega^0\times\Pc_{_2}(\R^d)\times\Ac$ $\mapsto$ $\rho_s^{t,\mu,\alpha}(\omega^0)$ 
$\in$  $\Pc_{_2}(\R^d)$ is measurable. 

 \vspace{1mm}

\noindent {\bf 3.} Let us finally check the flow property \reff{flowrho}. From pathwise uniqueness of the solution $\{X_s(\omega^0,.),t\leq s\leq T\}$  to \reff{Mckean} on $(\Omega,\Fc^1,\P^1)$ for $\P^0$-a.s. $\omega^0$ $\in$ 
$\Omega^0$, and recalling the definition of the shifted control process,  we have the flow property: for $t$ $\in$ $[0,T]$, $\xi$ $\in$ $L^2(\Fc_t.\R^d)$, $\alpha$ $\in$ $\Ac$, and $\P^0$-a.s. $\omega^0$ $\in$ 
$\Omega^0$,
\beqs
X_s^{t,\xi,\alpha}(\omega^0,.) &=& X_s^{\theta(\omega^0),X_{\theta(\omega^0)}^{t,\xi,\alpha}(\omega^0,.),\alpha^{\theta(\omega^0),\omega^0}}(\omega^0,.),  \;\;\; \P^1-\mbox{a.s.}
\enqs
for all $\F^0$-stopping time $\theta$ valued in $[t,T]$. It follows that for any Borel-measurable bounded function $\varphi$ on $\R^d$, and for  $\P^0$-a.s $\omega^0$ $\in$ $\Omega^0$,
\beqs
\rho_s^{t,\mu,\alpha}(\omega^0)(\varphi) \; = \;  \E^1 \Big[ \varphi\big(X_s^{t,\xi,\alpha}(\omega^0,.) \big) \Big]  & = & 
 \E^1 \Big[ \varphi\big(X_s^{\theta(\omega^0),X_{\theta(\omega^0)}^{t,\xi,\alpha}(\omega^0,.),\alpha^{\theta(\omega^0),\omega^0}}(\omega^0,.)\big) \Big]  \\
& = & \rho_s^{\theta(\omega^0),\rho_{\theta(\omega^0)}^{t,\mu,\alpha}(\omega^0),\alpha^{\theta(\omega^0),\omega^0}}(\omega^0)(\varphi), 
\enqs
where the last equality is obtained by noting  that $\rho_{\theta(\omega^0)}^{t,\mu,\alpha}(\omega^0)$ $=$ 
$\P_{X_{\theta(\omega^0)}^{t,\xi,\alpha}(\omega^0,.)}^{W^0}$, and the definition of $\rho_s^{t,\mu,\alpha}$.  This shows the required flow property \reff{flowrho}. 
\ep

\vspace{2mm}

%POINT: 
%--------  
% We need to prove the measurability of the FLOW  $(t,\mu,\alpha)$ $\mapsto$ $\rho^{t,\mu,\alpha}$ 
 %FOR this, we  have to prove in particular continuity with respect to $\alpha$ $\in$ $\Ac$ of $X^{t,\xi,\alpha}$.   
%---------  

Now, by the law of iterated conditional expectations, from \reff{defPW0}, \reff{defrho}, and recalling that $\alpha$ $\in$ $\Ac$ is $\F^0$-progressive, 
we can rewrite the cost functional as
\beq 
J(t,\xi,\alpha) &=& \E \Big[  \int_t^T \E \big[  f\big(X_s^{t,\xi},\P_{X_s^{t,\xi}}^{W^0},\alpha_s  \big) \big| \Fc_s^0 \big] ds 
+ \E \big[ g\big(X_T^{t,\xi},\P_{X_T^{t,\xi}}^{W^0}\big) \big| \Fc_T^0 \big]  \Big] \nonumber \\
&=& \E \Big[ \int_t^T \rho_s^{t,\mu} \big( f(.,\rho_s^{t,\mu},\alpha_s)\big) ds + \rho_T^{t,\mu}\big(g(.,\rho_T^{t,\mu})\big) \Big] \nonumber \\
&=& \E \Big[ \int_t^T \hat f(\rho_s^{t,\mu},\alpha_s) ds + \hat g(\rho_T^{t,\mu}) \Big], \label{Jmu}
\enq
for $t$ $\in$ $[0,T]$, $\xi$ $\in$ $L^2(\Gc;\R^d)$ with law $\mu$ $=$ $\Lc(\xi)$ $=$ $\P_\xi^{W^0}$  $\in$ $\Pc_{_2}(\R^d)$, $\alpha$ $\in$ $\Ac$, and with the functions $\hat f$ $:$ $\Pc_{_2}(\R^d)\times \bA$ $\rightarrow$ $\R$, and $\hat g$ $:$ $\Pc_{_2}(\R^d)$ $\rightarrow$ $\R$, defined by
\begin{equation}\label{hatfg}
\left\{
\begin{array}{ccc}
\hat f(\mu,a) &:=& \mu\big( f(.,\mu,a)\big) \; = \; \int_{\R^d} f(x,\mu,a) \mu(dx) \\
\hat g(\mu) &:=& \mu\big( g(.,\mu) \big) \; = \; \int_{\R^d} g(x,\mu) \mu(dx).
\end{array}
\right. 
\end{equation}
(To alleviate notations, we have omitted here the dependence of $\rho_s^{t,\mu}$ $=$ $\rho_s^{t,\mu,\alpha}$ on $\alpha$).  
Relation \reff{Jmu} means that the cost functional depends on $\xi$ only through its distribution $\mu$ $=$ $\Lc(\xi)$, and by misuse of notation, we set:
\beqs
J(t,\mu,\alpha)  \; := \;  J(t,\xi,\alpha) 
%& = & \E \Big[ \int_t^T \hat f(\rho_s^{t,\mu},\alpha_s) ds + \hat g(\rho_T^{t,\mu}) \Big], \\
&=& \E^0 \Big[ \int_t^T \hat f(\rho_s^{t,\mu},\alpha_s) ds + \hat g(\rho_T^{t,\mu}) \Big],  
\enqs
for $(t,\mu)$ $\in$ $[0,T]\times\Pc_{_2}(\R^d)$, $\xi \in L^2(\Gc;\R^d)$ with $\Lc(\xi)$ $=$ $\mu$, and the  expectation is taken under $\P^0$ since  
$\{\rho_s^{t,\mu},t\leq s\leq T\}$ is $\F^0$-progressive, and the control $\alpha$ $\in$ $\Ac$  is   an $\F^0$-progressive process. 
Therefore,  the value function can be  identified with a function defined on $[0,T]\times\Pc_{_2}(\R^d)$,  equal to (we keep the same notation $v(t,\mu)$ $=$ $v(t,\xi)$):
\beqs
v(t,\mu)  &=&  \inf_{\alpha\in\Ac} \E^0\Big[ \int_t^T \hat f(\rho_s^{t,\mu},\alpha_s) ds + \hat g(\rho_T^{t,\mu}) \Big],
\enqs
and satisfying from \reff{vquadra} the quadratic growth condition
\beq \label{vquadra2}
|v(t,\mu)| & \leq & C(1 + \|\mu\|_{_2}^2), \;\;\; \forall \mu \in \Pc_{_2}(\R^d). 
\enq

%by recalling that $\{\rho_s^{t,\mu},t\leq s\leq T\}$ is $\F^0$-progressive, and the control $\alpha$ $\in$ $\Ac$  is a random field $\F^0$-progressive process.   

\vspace{2mm}

As a consequence of the flow property in Lemma \ref{lemrho}, we obtain the following con\-ditioning lemma, also called pseudo-Markov property in the terminology of \cite{clataltan15}, for the 
controlled conditional distribution $\F^0$-progressive process $\{\rho_s^{t,\mu,\alpha},t\leq s\leq T\}$. 

\begin{Lemma} \label{lempseudo}
For any $(t,\mu,\alpha)$ $\in$ $[0,T]\times\Pc_{_2}(\R^d)\times\Ac$, and $\theta$ $\in$ $\Tc_{t,T}^0$, we have
\beq \label{Jcond}
%J(\theta(\omega^0),\rho_{\theta(\omega^0)}^{t,\mu,\alpha}(\omega^0),\alpha^{\theta(\omega^0),\omega^0}) &=&  
%\E^0 \Big[ \int_\theta^T \hat f(\rho_s^{t,\mu,\alpha},\alpha_s) ds + \hat g(\rho_T^{t,\mu,\alpha}) \big| \Fc_\theta^0 \Big](\omega^0), 
J(\theta,\rho_{\theta}^{t,\mu,\alpha},\alpha^{\theta}) &=&  
\E^0 \Big[ \int_\theta^T \hat f(\rho_s^{t,\mu,\alpha},\alpha_s) ds + \hat g(\rho_T^{t,\mu,\alpha}) \big| \Fc_\theta^0 \Big], \;\;\; \P^0-\mbox{a.s} 
\enq
%In particular, 
%\beq \label{vind}
%v(t,\mu)  &=&  \inf_{\alpha\in\Ac_t} \E^0\Big[ \int_t^T \hat f(\rho_s^{t,\mu},\alpha_s) ds + \hat g(\rho_T^{t,\mu}) \Big].
%\enq
\end{Lemma}
{\bf Proof.}  By the joint measurability property of $\rho_s^{t,\mu,\alpha}$  in $(t,s,\omega^0,\mu,\alpha)$ in Lemma \ref{lemrho},  the flow property \reff{flowrho}, and since 
$\rho_\theta^{t,\mu,\alpha}$ is $\Fc_\theta^0$-measurable for $\theta$ $\F^0$-stopping time,  we have for $\P^0$-a.s $\omega^0$ $\in$ $\Omega^0$, 
\beqs
& & \E^0 \Big[ \int_\theta^T \hat f(\rho_s^{t,\mu,\alpha},\alpha_s) ds + \hat g(\rho_T^{t,\mu,\alpha}) \big| \Fc_\theta^0 \Big](\omega^0) \\
&=&  \left. \E^0 \Big[ \int_r^T \hat f(\rho_s^{r,\pi,\beta},\beta_s) + \hat g(\rho_T^{r,\pi,\beta})  \big| \Fc_r^0 \Big](\omega^0)  
\right|_{r=\theta(\omega^0),\pi=\rho_{\theta(\omega^0)}^{t,\mu,\alpha}(\omega^0),\beta=\alpha^{r,\omega^0}} \\
&=& \left. \E^0 \Big[ \int_r^T \hat f(\rho_s^{r,\pi,\beta},\beta_s) + \hat g(\rho_T^{r,\pi,\beta})  \Big]  
\right|_{r=\theta(\omega^0),\pi=\rho_{\theta(\omega^0)}^{t,\mu,\alpha}(\omega^0),\beta=\alpha^{r,\omega^0}},
\enqs
where we used in the second equality the fact that for fixed $\omega^0$, $r$ $\in$ $[t,T]$, $\pi$ $\in$ $\Pc_{_2}(\R^d)$ represented by 
$\eta$ $\in$ $L^2(\Gc;\R^d)$ s.t. $\Lc(\xi)$ $=$ $\pi$, the process $\alpha^{r,\omega^0}$ lies in $\Ac_r$, hence is independent of $\Fc_r^0$, which implies that $X_s^{r,\eta,\alpha^{r,\omega^0}}$ is independent of $\Fc_r$, and thus $\rho_s^{r,\pi,\alpha^{r,\omega^0}}$ is also  
independent of $\Fc_r^0$ for $r\leq s$.  This shows the conditioning relation \reff{Jcond}. 
%Let us denote by $\tilde v(t,\mu)$ the r.h.s. of \reff{vind}. Since $\Ac_t$ $\subset$ $\Ac$, it is clear that $v(t,\mu)$ $\leq$ $\tilde v(t,\mu)$. To prove the reverse inequality, we apply the conditioning relation \reff{Jcond} for $\theta$ $=$ $t$, and get in %particular for all $\alpha$ $\in$ $\Ac$: 
%\beq \label{relvtildev}
%\int_{\Omega^0}   J(t,\mu,\alpha^{t,\omega^0})  \P^0(d\omega^0)  &=& J(t,\mu,\alpha). 
%\enq
%Now, recalling that for any fixed $\omega^0$ $\in$ $\Omega^0$, $\alpha^{t,\omega^0}$ lies in $\Ac_t$, we have $J(t,\mu,\alpha^{t,\omega^0})$ $\geq$ 
%$\tilde v(t,\mu)$, which proves the required result since $\alpha$ is arbitrary in \reff{relvtildev}.  
\ep

\subsection{Continuity of the value function and dynamic programming principle}

In this paragraph,  we show the continuity of the value function, which is helpful for proving next the dynamic programming principle. We mainly follow arguments from 
\cite{kry80} for the continuity result  that we extend  to our McKean-Vlasov framework.

\begin{Lemma} \label{contivJ}
The function $(t,\mu)$ $\mapsto$ $J(t,\mu,\alpha)$ is continuous on $[0,T]\times\Pc_{_2}(\R^d)$, uniformly with respect to  $\alpha$ $\in$ $\Ac$, and the function $\alpha$ $\mapsto$ $J(t,x,\alpha)$ is continuous on 
$\Ac$ for any $(t,\mu)$ $\in$ $[0,T]\times\Pc_{_2}(\R^d)$.  Consequently,  the cost functional $J$ is continuous on $[0,T]\times\Pc_{_2}(\R^d)\times\Ac$, and the value function $v$ is  
continuous on $[0,T]\times\Pc_{_2}(\R^d)$. 
\end{Lemma}
{\bf Proof.} {\it (1)}   For any  $0\leq t\leq s\leq T$, $\mu,\pi$ $\in$ $\Pc_{_2}(\R^d)$, $\alpha$ $\in$ $\Ac$,  recall that  $\P^0$-a.s. $\omega^0$ $\in$ $\Omega^0$, we have 
$\P^1_{X_r^{t,\xi,\alpha}(\omega^0,.)}$ $=$ $\rho_r^{t,\mu,\alpha}(\omega^0)$,  $\P^1_{X_r^{s,\zeta,\alpha}(\omega^0,.)}$ $=$ $\rho_r^{s,\pi,\alpha}(\omega^0)$ for $r$ $\in$ $[s,T]$, and any $\xi,\zeta$ $\in$ $L^2(\Gc;\R^d)$ s.t. 
$\Lc(\xi)$ $=$ $\mu$, $\Lc(\zeta)$ $=$ $\pi$.  By definition of $\|.\|_{_2}$ and  the Wasserstein distance  in $\Pc_{_2}(\R^d)$, we then have: $\|\rho_r^{t,\mu,\alpha}(\omega^0)\|_{_2}$ $=$ 
$\E^1|X_r^{t,\xi,\alpha}(\omega^0,.)|^2$,  and $\Wc_2^2\big(\rho_r^{t,\mu,\alpha}(\omega^0),\rho_r^{s,\pi,\alpha}(\omega^0)\big)$ $\leq$ 
$\E^1|X_r^{t,\xi,\alpha}(\omega^0,.)-X_r^{s,\zeta,\alpha}(\omega^0,.)|^2$, so that 
\beq 
\E^0 \Big[ \sup_{s\leq r\leq T} \|\rho_r^{t,\mu,\alpha}\|_{_2}^2  \Big] & \leq & \E \Big[ \sup_{s \leq  r \leq T} |X_r^{t, \xi, \alpha}|^2 \Big],  \label{estimXW0}  \\
\E^0 \Big[ \sup_{s\leq r\leq T} \Wc_2^2(\rho_r^{t,\mu,\alpha},\rho_r^{s,\pi,\alpha}) \Big] & \leq & \E\Big[\sup_{s \leq  r \leq T} |X_r^{t, \xi, \alpha}-X_r^{s, \zeta, \alpha}|^2 \Big].  \label{estimXW}
\enq
From the state equation \reff{Mckean}, and using standard arguments involving Burkholder-Davis-Gundy inequalities, \reff{estimXW0}, \reff{estimXW}, and Gronwall lemma, under the Lipschitz condition in {\bf(H1)}(i), 
we obtain the following estimates similar to the ones  for controlled diffusion processes (see \cite{kry80}, Chap.2, Thm.5.9, Cor.5.10):  there exists some positive constant $C$ s.t. for all $t$ 
$\in$ $[0,T]$,   $\xi, \zeta$ $\in$  $L^2(\Gc; \R^d)$, $\alpha$  $\in$ $\Ac$,  $h \in [0, T-t]$,
\beqs
%\E\big[\sup_{t \leq s \leq T} |X_s^{t, \xi, \alpha}|\big] &  \leq & C (1 +\E |\xi|^2),\\
\E\big[\sup_{t \leq s \leq t+h} |X_s^{t, \xi, \alpha}-\xi|^2\big] &  \leq & C(1+\E|\xi|^2)h,\\
\E\big[\sup_{t \leq s \leq T}|X_s^{t, \xi, \alpha}-X_s^{t, \zeta, \alpha}|^2\big] &  \leq&  C\E[|\xi -\zeta|^2],
\enqs
from which we easily deduce that for all $0\leq t\leq s\leq T$, $\xi, \zeta$ $\in$  $L^2(\Gc; \R^d)$, $\alpha$  $\in$ $\Ac$
\beq \label{estimX1}
\E\Big[\sup_{s \leq  r \leq T} |X_r^{t, \xi, \alpha}-X_r^{s, \zeta, \alpha}|^2 \Big] &  \leq &  C \big(\E|\xi-\zeta|^2 +(1 + \E|\xi|^2 + \E|\zeta|^2)|s-t| \big). 
\enq
Together with the estimates \reff{momentX}, and by definition of $\Wc_{_2}(\mu,\pi)$, $\|\mu\|_{_2}$, $\|\pi\|_{_2}$, we then get from \reff{estimXW0}, \reff{estimXW}: 
\beq 
\E^0 \Big[ \sup_{s\leq r\leq T} \|\rho_r^{t,\mu,\alpha}\|_{_2}^2  \Big] & \leq & C(1 + \|\mu\|^2_{_2}),   \label{estimrho0} \\
\E^0 \Big[ \sup_{s\leq r\leq T} \Wc_2^2(\rho_r^{t,\mu,\alpha},\rho_r^{s,\pi,\alpha}) \Big] & \leq & C \big( \Wc_2^2(\mu,\pi) + (1 + \|\mu\|_{_2}^2 +  \|\pi\|_{_2}^2)|s-t| \big).   \label{estimrho1}
\enq
{\it (2)}  Let us now show the continuity of  the cost functional  $J$ in $(t,\mu)$, uniformly w.r.t. $\alpha$ $\in$ $\Ac$. First, we notice  from the growth condition in {\bf (H2)}(i) and the local Lipschitz condition in {\bf (H2)}(ii) 
that there exists some positive constant $C$ s.t. 
for all $\mu,\pi$ $\in$ $\Pc_{_2}(\R^d)$, $\alpha$ $\in$ $\Ac$, 
\beqs
|\hat f(\mu,\alpha)| & \leq & C\big(1 + \|\mu\|_{_2}^2\big),  \\
|\hat f(\mu,\alpha) - \hat f(\pi,\alpha)| + |\hat g(\mu)-\hat g(\pi)| & \leq & C(1 + \|\mu\|_{_2} + \|\pi\|_{_2})  \Wc_2(\mu,\pi). 
\enqs
Then, we have for all $0\leq t\leq s\leq T$, $\mu,\pi$ $\in$ $\Pc_{_2}(\R^d)$, $\alpha$ $\in$ $\Ac$
\beqs
\big| J(t,\mu,\alpha) - J(s,\pi,\alpha) \big| & \leq & \E^0 \Big[ \int_t^s  |\hat f(\rho_r^{t,\mu,\alpha})| dr \Big] \\
& &   + \; \E^0 \Big[ \int_s^T \big| \hat f(\rho_r^{t,\mu,\alpha},\alpha_r) -  \hat f(\rho_r^{s,\pi,\alpha},\alpha_r) \big| dr 
+ \big| \hat g(\rho_T^{t,\mu,\alpha}) -  \hat g(\rho_T^{s,\pi,\alpha}) \big| \Big] \\
& \leq & C\E^0\Big[ (1 +  \sup_{t\leq r\leq s} (\|\rho_r^{t,\mu,\alpha}\|_{_2})|s-t| \Big] \\ 
& &  + \; C   \E^0 \Big[ \big(1 +  \sup_{s\leq r\leq T} (\|\rho_r^{t,\mu,\alpha}\|_{_2}  + \|\rho_r^{s,\pi,\alpha}\|_{_2})  \big) \sup_{s\leq r\leq T} \Wc_2(\rho_r^{t,\mu,\alpha},\rho_r^{s,\pi,\alpha}) \Big]  \\
& \leq &   C (1 + \|\mu\|_{_2}) |s-t| \\
& &  + \; C (1 + \|\mu\|_{_2} +  \|\pi\|_{_2}) \big( \Wc_2(\mu,\pi) + (1 + \|\mu\|_{_2} +  \|\pi\|_{_2})|s-t|^{1\over 2} \big),
\enqs
by Cauchy Schwarz inequality  and \reff{estimrho0}-\reff{estimrho1}, which shows the desired  continuity result.

 \vspace{1mm}

 \noindent {\it (3)}  Let us show the continuity of the cost functional with respect to the control. Fix  $(t,\mu)$ $\in$ $[0,T]\times\Pc_{_2}(\R^d)$, and consider   $\alpha$ $\in$ $\Ac$, 
 a sequence  $(\alpha^n)_n$ in $\Ac$ s.t. $\Delta(\alpha^n,\alpha)$ $\rightarrow$ $0$, i.e.  $d_A(\alpha^n_t,\alpha_t)$ $\rightarrow$ $0$ in $dt\otimes d\P^0$-measure, as $n$ goes to infinity.   
 Denote by  $\rho^n$ $=$ $\rho^{t,\mu,\alpha^n}$, $\rho$ $=$ $\rho^{t,\mu,\alpha}$, $X^n$ $=$ $X^{t,\xi,\alpha^n}$, $X$ $=$ $X^{t,\xi,\alpha}$ for  $\xi$ $\in$ $L^2(\Gc;\R^d)$ s.t. $\Lc(\xi)$ $=$ $\mu$. 
 By the same arguments as in \reff{estimXW}, we have
 \beq \label{estimXWn}
 \E^0 \Big[ \sup_{t\leq s\leq T} \Wc_2^2(\rho_s^n,\rho_s) \Big] & \leq & \E\Big[\sup_{t \leq  s \leq T} |X_s^{n}-X_s|^2 \Big]. 
 \enq
 Next, starting from the state equation \reff{Mckean}, using standard arguments involving Burkholder-Davis-Gundy inequalities, \reff{estimXWn}, and Gronwall lemma, under the Lipschitz condition in {\bf(H1)}(i), we arrive at:
 \beqs
 \E\Big[\sup_{t \leq  s \leq T} |X_s^{n}-X_s|^2 \Big] & \leq & C \Big\{ \E\Big[ \int_t^T |b(X_s,\rho_s,\alpha_s)- b(X_s,\rho_s,\alpha_s^n)|^2  ds \\
 & & \;\;\;\;\;\;\;  + \;  \int_t^T |\sigma(X_s,\rho_s,\alpha_s)- \sigma(X_s,\rho_s,\alpha_s^n)|^2  ds \\
 & &  \;\;\;\;\;\;\;  + \;  \int_t^T |\sigma_0(X_s,\rho_s,\alpha_s)- \sigma_0(X_s,\rho_s,\alpha_s^n)|^2  ds \Big] \Big\},
 \enqs
 for some positive constant $C$ independent of $n$. Recalling the bound \reff{momentX}, and \reff{estimXW0}, we deduce by the dominated convergence theorem under  the  linear growth condition in {\bf (H1)}(i), and the continuity assumption in {\bf (H1)}(ii) that   $\E\big[\sup_{t \leq  s \leq T} |X_s^{n}-X_s|^2 \big]$ $\rightarrow$ $0$,  and thus by \reff{estimXWn}
\beq \label{convrhon} 
\E^0 \Big[ \sup_{t\leq s\leq T} \Wc_2^2(\rho_s^n,\rho_s) \Big] & \rightarrow & 0, \;\;\; \mbox{ as } \; n \rightarrow \infty. 
\enq 
Now, by writing
\beq 
& & \big| J(t,\mu,\alpha^n) - J(t,\mu,\alpha) \big| \nonumber \\
& \leq &  \E^0 \Big[ \int_t^T \big| \hat f(\rho_s^{n},\alpha_s^n) -  \hat f(\rho_s,\alpha_s) \big| ds 
+ \big| \hat g(\rho_T^{n}) -  \hat g(\rho_T) \big| \Big], \label{Jcont}
\enq
and noting that $\hat f$ and $\hat g$ are continuous on $\Pc_{_2}(\R^d)\times\bA$, resp. on $\Pc_{_2}(\R^d)$,  under the continuity assumption in {\bf (H2)}(ii),  we conclude by the same arguments as in 
\cite{kry80} using \reff{convrhon} (see Chapter 3, Sec. 2,  or  also  Lemma 4.1 in \cite{fuhpha15}) that   the r.h.s. of \reff{Jcont} tends to zero as $n$ goes to infinity, which proves the continuity of $J(t,\mu,.)$ on $\Ac$.

\vspace{1mm}

\noindent {\it (4)}  Finally, the global continuity of the cost functional $J$ on $[0,T]\times\Pc_{_2}(\R^d)\times\Ac$ is a direct consequence of the continuity of $J(.,.,\alpha)$ on $[0,T]\times\Pc_{_2}(\R^d)$ 
uniformly w.r.t. $\alpha$ $\in$ $\Ac$, and the continuity of $J(t,\mu,.)$ on $\Ac$, while the continuity of the value function $v$ on $[0,T]\times\Pc_{_2}(\R^d)$ follows immediately from the fact that 
\beqs
|v(t,\mu) - v(s,\pi)| & \leq & \sup_{\alpha\in\Ac} |J(t,\mu,\alpha) - J(s,\pi,\alpha) |, \;\;\;  t,s \in [0,T], \; \mu,\pi \in \Pc_{_2}(\R^d), 
\enqs
 and again from  the continuity of $J(.,.,\alpha)$ on $[0,T]\times\Pc_{_2}(\R^d)$ uniformly w.r.t. $\alpha$ $\in$ $\Ac$.  
 \ep 
 
\vspace{1mm}

\begin{Remark}
{\rm Notice that the supremum defining the value function $v(t,\mu)$ can be taken over the subset $\Ac_t$ of  elements in $\Ac$ which are independent of $\Fc_t^0$ under $\P^0$, i.e.
\beq \label{vind}
v(t,\mu)  &=&  \inf_{\alpha\in\Ac_t} \E^0\Big[ \int_t^T \hat f(\rho_s^{t,\mu},\alpha_s) ds + \hat g(\rho_T^{t,\mu}) \Big].
\enq
Indeed, denoting by $\tilde v(t,\mu)$ the r.h.s. of \reff{vind}, and since $\Ac_t$ $\subset$ $\Ac$, it is clear that $v(t,\mu)$ $\leq$ $\tilde v(t,\mu)$. To prove the reverse inequality, we apply the conditioning relation 
\reff{Jcond} for $\theta$ $=$ $t$, and get in particular for all $\alpha$ $\in$ $\Ac$: 
\beq \label{relvtildev}
\int_{\Omega^0}   J(t,\mu,\alpha^{t,\omega^0})  \P^0(d\omega^0)  &=& J(t,\mu,\alpha). 
\enq
Now, recalling that for any fixed $\omega^0$ $\in$ $\Omega^0$, $\alpha^{t,\omega^0}$ lies in $\Ac_t$, we have $J(t,\mu,\alpha^{t,\omega^0})$ $\geq$ 
$\tilde v(t,\mu)$, which proves the required result since $\alpha$ is arbitrary in \reff{relvtildev}.  
}
\ep
\end{Remark}

\vspace{2mm}

We can now state the dynamic programming principle (DPP) for the value function to the stochastic McKean-Vlasov control problem.

\begin{Proposition} \label{propdyn} (Dynamic Programming Principle)

\noindent We have for all  $(t,\mu)$ $\in$ $[0,T]\times\Pc_{_2}(\R^d)$, 
%and $\theta$ $\in$ $\Tc^0_{t,T}$:
%\beq \label{DPPdeter}
%v(t,\mu) &=& \inf_{\alpha\in\Ac} \E^0\Big[ \int_t^\theta \hat f(\rho_s^{t,\mu,\alpha},\alpha_s)  ds \; + \; v(\theta,\rho_\theta^{t,\mu,\alpha}) \Big]. 
%\enq
\beqs 
v(t,\mu) &=& \inf_{\alpha\in\Ac} \inf_{\theta\in\Tc_{t,T}^0} 
\E^0\Big[ \int_t^\theta \hat f(\rho_s^{t,\mu,\alpha},\alpha_s)  ds \; + \; v(\theta,\rho_\theta^{t,\mu,\alpha}) \Big]
\label{DPPstrong1} \\
&=& \inf_{\alpha\in\Ac} \sup_{\theta\in\Tc_{t,T}^0} 
\E^0\Big[ \int_t^\theta \hat f(\rho_s^{t,\mu,\alpha},\alpha_s)  ds \; + \; v(\theta,\rho_\theta^{t,\mu,\alpha}) \Big],
\enqs
which means equivalently that

\noindent (i) for all $\alpha$ $\in$ $\Ac$, $\theta$ $\in$ $\Tc_{t,T}^0$, 
\beq \label{DPP1}
v(t,\mu) & \leq & \E^0\Big[ \int_t^\theta \hat f(\rho_s^{t,\mu,\alpha},\alpha_s)  ds \; + \; v(\theta,\rho_\theta^{t,\mu,\alpha}) \Big],
\enq
(ii) for all $\eps$ $>$ $0$, there exists $\alpha$ $\in$ $\Ac$, such that for all $\theta$ $\in$ $\Tc_{t,T}^0$, 
\beq \label{DPP2}
v(t,\mu) + \eps & \geq & \E^0\Big[ \int_t^\theta \hat f(\rho_s^{t,\mu,\alpha},\alpha_s)  ds \; + \; v(\theta,\rho_\theta^{t,\mu,\alpha}) \Big]. 
\enq
\end{Proposition}

\begin{Remark}
{\rm The above formulation of the DPP implies in particular that for all $\theta$ $\in$ $\Tc_{t,T}^0$, 
\beqs \label{dynproweak}
v(t,\mu) & =  &  \inf_{\alpha\in\Ac} \E^0\Big[ \int_t^\theta \hat f(\rho_s^{t,\mu,\alpha},\alpha_s)  ds \; + \; v(\theta,\rho_\theta^{t,\mu,\alpha}) \Big],
\enqs
which is the usual formulation of the DPP.  The formulation in Proposition \ref{propdyn} is stronger, and the difference relies on the fact that in the inequality \reff{DPP2}, the $\eps$-optimal control $\alpha$ $=$ 
$\alpha^\eps$ does not depend on $\theta$. This condition will be useful to show later the viscosity supersolution property of the value function. 
\ep
}
\end{Remark}

\noindent {\bf Proof.} {\bf 1.} Fix  $(t,\mu)$ $\in$ $[0,T]\times\Pc_{_2}(\R^d)$. From the conditioning relation \reff{Jcond}, we have 
for all $\theta$ $\in$ $\Tc_{t,T}^0$, $\alpha$ $\in$ $\Ac$, 
\beq \label{Jtheta}
J(t,\mu,\alpha) &=& \E^0\Big[ \int_t^\theta \hat f(\rho_s^{t,\mu,\alpha},\alpha_s)  ds + J(\theta,\rho_\theta^{t,\mu,\alpha},\alpha^\theta) \Big]. 
\enq
Since $J(.,.,\alpha^\theta)$ $\geq$ $v(.,.)$, and $\theta$ is arbitrary in $\Tc^0_{t,T}$, we have
\beqs
J(t,\mu,\alpha) & \geq & \sup_{\theta\in\Tc_{t,T}^0} 
\E^0\Big[ \int_t^\theta \hat f(\rho_s^{t,\mu,\alpha},\alpha_s)  ds + v(\theta,\rho_\theta^{t,\mu,\alpha}) \Big],
\enqs
and since  $\alpha$ is arbitrary in $\Ac$, it follows that
\beq 
v(t,\mu) & \geq & \inf_{\alpha\in\Ac} \sup_{\theta\in\Tc_{t,T}^0} 
\E^0\Big[ \int_t^\theta \hat f(\rho_s^{t,\mu,\alpha},\alpha_s)  ds \; + \; v(\theta,\rho_\theta^{t,\mu,\alpha}) \Big]
\label{DPPstrong1inter}
\enq

\noindent {\bf 2.}  Fix $(t,\mu)$ $\in$ $[0,T]\times\Pc_{_2}(\R^d)$, $\alpha$ $\in$ $\Ac$ and $\theta$ $\in$ $\Tc_{t,T}^0$.  For any $\eps$ $>$ $0$, 
$\omega^0$ $\in$ $\Omega^0$, one can find from \reff{vind} some $\alpha^{(\eps,\omega^0)}$ $\in$ $\Ac_{\theta(\omega^0)}$ s.t.
\beq \label{vthetaJ}
v(\theta(\omega^0),\rho_{\theta(\omega^0)}^{t,\mu,\alpha}(\omega^0)) + \eps & \geq & 
J(\theta(\omega^0),\rho_{\theta(\omega^0)}^{t,\mu,\alpha}(\omega^0),\alpha^{(\eps,\omega^0)}). 
\enq
Since $J$ and  $v$ are continuous (by Lemma \ref{contivJ}), one can invoke measurable selection arguments (see e.g. \cite{wag80}), to claim  that the map $\omega^0$ $\in$ $(\Omega^0,\Fc^0)$ 
$\mapsto$ $\alpha^{(\eps,\omega^0)}$ $\in$ $(\Ac,\Bc_\Ac)$ can be chosen measurable.  
Let us now define the process $\bar\alpha$ on $(\Omega^0,\Fc^0,\P^0)$ 
obtained by concatenation at $\theta$ of the processes $\alpha$ and $\alpha^{(\eps,\omega^0)}$ in $\Ac$, namely:
\beqs
\bar\alpha_s(\omega^0) &:=& \alpha_s(\omega^0) 1_{s<\theta(\omega^0)} + \alpha^{(\eps,\omega^0)}(\omega^0) 1_{s\geq\theta(\omega^0)}, \;\;\; 
0 \leq s \leq T.
\enqs
By Lemma 2.1 in \cite{sontou02}, and since $\Ac$ is a separable metric space,  the process $\bar\alpha$ is  $\F^0$-progressively measurable, 
and thus $\bar\alpha$  $\in$ $\Ac$.  Notice with our notations of shifted control process that $\bar\alpha^{\theta(\omega^0),\omega^0}$ $=$ 
$\alpha^{(\eps,\omega^0)}$ for all $\omega^0$ in $\Omega^0$, and then \reff{vthetaJ} reads as
\beqs
v(\theta,\rho_\theta^{t,\mu,\alpha}) + \eps & \geq & J(\theta,\rho_\theta^{t,\mu,\alpha},\bar\alpha^\theta), \;\;\; \P^0-\mbox{a.s.}
\enqs
Therefore, by using again \reff{Jtheta} to $\bar\alpha$, and since $\rho_s^{t,\mu,\bar\alpha}$ $=$ $\rho_s^{t,\mu,\alpha}$ for $s$ $\leq$ $\theta$ (recall that 
$\bar\alpha_s$ $=$ $\alpha_s$ for $s$ $<$ $\theta$, and $\rho^{t,\mu,\alpha}$ has continuous trajectories),  we get
\beqs
v(t,\mu)  \leq \;  J(t,\mu,\bar\alpha) & =  & 
\E^0\Big[ \int_t^\theta \hat f(\rho_s^{t,\mu,\alpha},\alpha_s)  ds + J(\theta,\rho_\theta^{t,\mu,\alpha},\bar\alpha^\theta) \Big] \\
& \leq &  \E^0\Big[ \int_t^\theta \hat f(\rho_s^{t,\mu,\alpha},\alpha_s)  ds + v(\theta,\rho_\theta^{t,\mu,\alpha}) \Big]  + \eps
\enqs
Since $\alpha$, $\theta$ and $\eps$ are arbitrary, this gives the inequality
\beqs
v(t,\mu) & \leq & \inf_{\alpha\in\Ac} \inf_{\theta\in\Tc_{t,T}^0} 
\E^0\Big[ \int_t^\theta \hat f(\rho_s^{t,\mu,\alpha},\alpha_s)  ds \; + \; v(\theta,\rho_\theta^{t,\mu,\alpha}) \Big],
\enqs
which, combined with the first inequality \reff{DPPstrong1inter}, proves the DPP result.
\ep

\section{Bellman equation and viscosity solutions}

\setcounter{equation}{0} \setcounter{Assumption}{0}
\setcounter{Theorem}{0} \setcounter{Proposition}{0}
\setcounter{Corollary}{0} \setcounter{Lemma}{0}
\setcounter{Definition}{0} \setcounter{Remark}{0}

\subsection{Differentiability and It\^o's formula in Wasserstein space}

We shall rely on the notion of derivative with respect to a probability measure, as introduced by P.L. Lions in his course at Coll\`ege de France \cite{lio12}. 
We provide a brief introduction to this concept  and refer to  the lecture notes  \cite{car12} (see also \cite{buetal14}, \cite{chacridel15}) for the details.

This notion is based on the {\it lifting}  of functions $u$ $:$ $\Pc_{_2}(\R^d)$ $\rightarrow$ $\R$ into functions $\tilde u$ defined on 
$L^2(\Gc;\R^d)$ ($=$ $L^2(\tilde\Omega^1,\Gc,\tilde\P^1;\R^d)$)  by setting $\tilde u(\xi)$ $=$ $u(\Lc(\xi))$ ($=$ $u(\tilde\P^1_{\xi})$). 
Conversely, given a function $\tilde u$ defined on $L^2(\Gc;\R^d)$,  we call {\it inverse-lifted} function of $\tilde u$ the function $u$ defined on $\Pc_{_2}(\R^d)$ by 
$u(\mu)$ $=$ $\tilde u(\xi)$ for $\mu$ $=$ $\Lc(\xi)$, and we notice that such $u$ exists iff $\tilde u(\xi)$ depends only on the distribution of $\xi$ for any $\xi$ $\in$ $L^2(\Gc;\R^d)$. In this case, we shall often identify in the sequel 
the function $u$ and its lifted version $\tilde u$, by using the same notation $u$ $=$ $\tilde u$.  

We say that $u$ is differentiable (resp. $\Cc^1$) on
$\Pc_{_2}(\R^d)$ if the lift $\tilde u$ is Fr\'echet differentiable (resp. Fr\'echet differentiable with continuous derivatives) on  $L^2(\Gc;\R^d)$. 
In this case, the Fr\'echet derivative $[D \tilde u](\xi)$, viewed as an element $D\tilde u(\xi)$ of $L^2(\Gc;\R^d)$  by Riesz' theorem:  
$[D \tilde u](\xi)(Y)$ $=$ $\tilde\E[D\tilde u(\xi).Y]$,  can be represented as
\beq \label{Uu1}
D\tilde u(\xi) &=& \partial_\mu u(\Lc(\xi))(\xi),
\enq
for some function  $\partial_\mu u(\Lc(\xi))$ $:$ $\R^d$ $\rightarrow$ $\R^d$,  which is  called derivative of $u$ at $\mu$ $=$ $\Lc(\xi)$.  Moreover, 
$\partial_\mu u(\mu)$ $\in$ $L^2_\mu(\R^d)$ for  $\mu$ $\in$ $\Pc_{_2}(\R^d)$ $=$ $\{\Lc(\xi):  \xi  \in L^2(\Gc;\R^d)\}$. 
Following \cite{chacridel15}, we say that $u$ is fully  $\Cc^2$ if it is $\Cc^1$, 
and one can find, for any $\mu$ $\in$ $\Pc_{_2}(\R^d)$, a continuous version of the mapping $x\in\R^d$ $\mapsto$ $\partial_\mu u(\mu)(x)$, such that the mapping
$(\mu,x)$ $\in$ $\Pc_{_2}(\R^d)\times\R^d$ $\mapsto$ $\partial_\mu u(\mu)(x)$  is continuous at any point $(\mu,x)$ such that $x$ $\in$ Supp$(\mu)$, and 
%the mapping $(\mu,x)$ $\in$ $\Pc_{_2}(\R^d)\times\R^d$ $\mapsto$ $\partial_\mu u(\mu)(x)$  is continuous and
\begin{itemize}
\item[(i)] for each fixed $\mu$ $\in$ $\Pc_{_2}(\R^d)$, the mapping $x$ $\in$ $\R^d$ $\mapsto$  $\partial_\mu u(\mu)(x)$ is differentiable in the standard sense, with a gradient denoted by  
$\partial_x  \partial_\mu u(\mu)(x)$  $\in$ $\R^{d\times d}$, and s.t. the mapping  $(\mu,x)$ $\in$ $\Pc_{_2}(\R^d)\times\R^d$ 
$\mapsto$  $\partial_x  \partial_\mu u(\mu)(x)$ is continuous
\item[(ii)] for each fixed  $x$ $\in$ $\R^d$, the mapping $\mu$ $\in$ $\Pc_{_2}(\R^d)$ $\mapsto$  $\partial_\mu u(\mu)(x)$ is differentiable in the above lifted sense.  Its derivative, interpreted thus as a mapping $x'$ $\in$ $\R^d$ $\mapsto$ $\partial_\mu \big[ \partial_\mu u(\mu)(x)\big](x')$ $\in$ $\R^{d\times d}$ in 
$L^2_\mu(\R^{d\times d})$, is denoted by $x'$ $\in$ $\R^d$ $\mapsto$ $\partial_\mu^2 u(\mu)(x,x')$, and s.t. the mapping $(\mu,x,x')$ $\in$ 
$\Pc_{_2}(\R^d)\times\R^d\times\R^d$ $\mapsto$ $\partial_\mu^2 u(\mu)(x,x')$ is continuous. 
\end{itemize}
We say that $u$ $\in$ $\Cc^2_b(\Pc_{_2}(\R^d))$ if it is fully $\Cc^2$,  $\partial_x  \partial_\mu u(\mu)$ $\in$ $L_\mu^\infty(\R^{d\times d})$, $\partial_\mu^2 u(\mu)$ $\in$ $L_{\mu\otimes\mu}^\infty(\R^{d\times d})$ for any $\mu$ 
$\in$ $\Pc_{_2}(\R^d)$,  and for any compact set $\Kc$ of $\Pc_{_2}(\R^d)$, we have
%there exists some constant $C_u$ s.t. for all $\mu$ $\in$ $\Pc_{_2}(\R^d)$,
\beq \label{Kbor}
 \sup_{ \mu \in \Kc } \Big[ \int_{\R^d} \big| \partial_\mu u(\mu)(x) |^2\mu(dx)  +
%\int_{\R^d}  \big| \partial_x \partial_\mu u(\mu)(x) |^2 \mu(dx)  
\big \| \partial_x \partial_\mu u(\mu)\|_{_\infty} +  \big \| \partial_\mu^2 u(\mu)\|_{_\infty}
\Big]  & < & \infty.
%& \leq & C_u\big( 1 + \|\mu\|^2_{_2} \big).
\enq
We next need an It\^o's formula along a flow of conditional measures proved in \cite{chacridel15} (see also \cite{caretal16} and \cite{cardel14}).  
%for processes with common noise. 
Let  $(\Omega,\Fc,\P)$ be a probability space of the form  $(\Omega,\Fc,\P)$ $=$ $(\Omega^{0}\times\Omega^1,\Fc^{0}\otimes\Fc^1,\P^0\otimes\P^1)$, where $(\Omega^0,\Fc^0,\P^0)$ supports $W^0$ and 
$(\Omega^1,\Fc^1,\P^1)$ supports $B$ as in Section \ref{seccondMcKean}.   
Let us consider an It\^o process in $\R^d$  of the form:
\beq \label{ItoX}
X_t &=& X_0 +  \int_0^t b_s ds +  \int_0^t \sigma_s dB_s +  \int_0^t \sigma_s^0 dW^0_s, \;\;\; 0 \leq t \leq T, 
\enq
where  $X_0$ is independent of $(B,W^0)$,  $b$, $\sigma$, $\sigma^0$ are progressively measurable processes with respect to the natural filtration $\F$ generated by 
$(X^0,B,W^0)$, and satisfying the square integrability condition: $\E \big[ \int_0^T |b_t|^2 + |\sigma_t|^2 + |\sigma_t^0|^2 dt \big]$ $<$ $\infty$.  
Denote by $\P_{_{X_t}}^{_{W^0}}$ the conditional law of $X_t$, $t$ $\in$ $[0,T]$, given the $\sigma$-algebra $\Fc^0$ generated by the 
whole filtration of $W^0$, and by $\E_{_{W^0}}$ $=$ $\E^1$ the conditional expectation w.r.t. $\Fc^{0}$. 
Let $u$ $\in$ $\Cc^2_b(\Pc_{_2}(\R^d))$. Then, for all $t$ $\in$ $[0,T]$, we have: 
\beq
u(\P_{_{X_t}}^{_{W^0}}) &=& u(\P_{_{X_0}})  + \int_0^t \E_{_{W^0}} \Big[ \partial_\mu u(\P_{_{X_s}}^{_{W^0}})(X_s).b_s 
+ \frac{1}{2}{\rm tr}\big(\partial_x \partial_\mu u(\P_{_{X_s}}^{_{W^0}})(X_s) (\sigma_s\sigma_s\trans + \sigma_s^0(\sigma_s^0)\trans) \big) \Big] \nonumber \\
& & \;\;\;\;\;\;\;  + \;  \E_{_{W^0}} \Big[ \E'_{_{W^0}} 
\big[ \frac{1}{2}{\rm tr}\big(\partial_\mu^2 u(\P_{_{X_s}}^{_{W^0}})(X_s,X_s') \sigma_s^0(\sigma_s^{'0})\trans\big) \big] \Big] ds  \nonumber \\
& & \;\;\; + \; \int_0^t \E_{_{W^0}} \Big[ \partial_\mu  u(\P_{_{X_s}}^{_{W^0}})(X_s)\trans \sigma_s^0 \Big] dW_s^0, \label{Ito}
\enq
where  $X'$ and $\sigma^{'0}$ are copies of  $X$ and $\sigma^0$ on another probability space 
$(\Omega'=\Omega^0\times\Omega^{'1}\,\Fc^0\otimes\Fc^{'1},\P^0\times\P^{'1})$, with  $(\Omega^{'1},\Fc^{'1},\P^{'1})$ supporting  $B'$ a copy of $B$, and $\E'_{_{W^0}}$ $=$ $\E^{'1}$.   Here $\trans$ denotes the transpose of any vector or matrix.

%In this case, notice that $\partial_\mu u(\mu)$ (resp.  $\partial_x \partial_\mu u(\mu)$)  lies in $L_\mu^2(\R^d;\R^d)$ (resp. $L_\mu^2(\R^d;\R^{d\times d})$) the set functions from $\R^d$ into $\R^d$ (resp. $\R^{d\times d}$),  
%and square-integrable w.r.t. $\mu$.  

In the sequel, it will be useful to formulate It\^o's formula  for the lifted function $\tilde u$ on 
$L^2(\Gc,\R^d)$ ($=$ $L^2(\tilde\Omega^1,\Gc,\tilde\P^1;\R^d)$). Notice, however, that even if $u$ $\in$ 
$\Cc^2_b(\Pc_{_2}(\R^d))$, then its  lifted function $\tilde u$  may not be in general  twice continuously Fr\'echet differentiable on   $L^2(\Gc;\R^d)$, as 
discussed in Example 2.1 in \cite{buetal14}. Under the extra-assumption that the lift $\tilde u$ $\in$ $\Cc^2(L^2(\Gc;\R^d))$,  the second Fr\'echet derivative $D^2 \tilde u(\xi)$  is identified indifferently by Riesz' theorem  as a  bilinear form on $L^2(\Gc;\R^d)$ or  as a self-adjoint operator  (hence bounded) on  $L^2(\Gc;\R^d)$, denoted  by  $D^2\tilde u(\xi)$ $\in$ $S(L^2(\Gc;\R^d))$,  and we have the relation (see Appendix A.2  in \cite{cardel14b}): 
\begin{equation} \label{Uu2}
\left\{
\begin{array}{ccccl}
 D^2\tilde u(\xi)[Y, Y] & = & \tilde\E^1 \Big[ D^2\tilde u(\xi)(Y).Y \Big]  &=& \tilde\E^1 \Big[ \tilde\E^{'1} \big[ {\rm tr} \big( \partial_\mu^2 u(\Lc(\xi))(\xi,\xi') Y(Y')\trans \big) \big] \Big] \\
&&  & & \;\;\;\;\;  + \;\;  \tilde\E^1\Big[ {\rm tr} \big( \partial_x \partial_\mu u(\Lc(\xi))(\xi) YY\trans \big) \Big], \\
 D^2\tilde u(\xi)[ZN, ZN]  & = & \tilde\E^1 \Big[ D^2\tilde u(\xi)(ZN).ZN \Big]  &=& \tilde\E^1\Big[ {\rm tr} \big( \partial_x \partial_\mu u(\Lc(\xi))(\xi) ZZ\trans \big) \Big], 
 \end{array}
 \right.
\end{equation}
for any $\xi$ $\in$ $L^2(\Gc;\R^d)$, $Y$ $\in$ $L^2(\Gc;\R^{d})$,  $Z$ $\in$ $L^2(\Gc;\R^{d\times q})$, and where $(\xi',Y')$ is a copy of  
$(\xi,Y)$ on another Polish and atomless probability space $(\tilde\Omega^{'1},\Gc',\tilde\P^{'1})$, $N$ $\in$ $L^2(\Gc;\R^q)$ is  independent of 
$(\xi,Z)$ with zero mean, and unit variance. Now,  let is  consider a copy $\tilde B$ of $B$ on 
the probability space $(\tilde\Omega^1,\Gc,\tilde\P^1)$, denote by $\tilde X_0$, $\tilde b$, $\tilde\sigma$, $\tilde\sigma_0$ copies of $X_0$, $b$, 
$\sigma$, $\sigma_0$ on 
$(\tilde\Omega=\Omega^0\times\tilde\Omega^1,\tilde\Fc=\Fc^0\otimes\Gc,\tilde\P =\P^0\otimes\tilde\P^1)$, and consider the It\^o process $\tilde X$ on $(\tilde\Omega,\tilde\Fc,\tilde P)$ of the form
\beqs
\tilde X_t &=& \tilde X_0 +  \int_0^t \tilde b_s ds +  \int_0^t \tilde \sigma_s d\tilde B_s +  \int_0^t \tilde \sigma_s^0 dW^0_s, \;\;\; 0 \leq t \leq T, 
\enqs
which is then a copy of $X$ in \reff{ItoX}.  The process $\check X$ defined by $\check X_t(\omega^0)$ $=$ $\tilde X_t(\omega^0,.)$, $0\leq t\leq T$,  is $\F^0$-progressive, and valued in 
$L^2(\Gc;\R^d)$.  Similarly, the processes defined by $\check b_t(\omega^0)$ $=$ $\tilde b_t(\omega^0,.)$, $\check \sigma_t(\omega^0)$ $=$ $\tilde \sigma_t(\omega^0,.)$, 
$\check \sigma^0_t(\omega^0)$ $=$ $\tilde \sigma^0_t(\omega^0,.)$, $0\leq t\leq T$, are valued in $L^2(\Gc;\R^d)$, $\P^0$-a.s.  Thus, when 
the lifted function $\tilde u$ $\in$ $\Cc^2(L^2(\Gc;\R^d))$, we obtain from \reff{Ito}  and relation \reff{Uu1}-\reff{Uu2} an It\^o's formula on the 
lifted space $L^2(\Gc;\R^d)$:
\beq
\tilde u(\check X_t) &=& \tilde u(\check X_0) +  \int_0^t \tilde\E^1\Big[ D\tilde u(\check X_s).\check b_s + \frac{1}{2} D^2\tilde u(\check X_s)(\check\sigma_s N).\check\sigma_s N + 
 \frac{1}{2} D^2\tilde u(\check X_s)(\check\sigma_s^0).\check\sigma_s^0   \Big] ds \nonumber \\
 & & \;\;\; + \; \int_0^t \tilde\E^1 \big[ D\tilde u(\check X_s)\trans\check\sigma_s^0 \big] dW_s^0, \;\;\;\;\;\;\;  0 \leq t \leq T,   \; \P^0-a.s. \label{Ito2}
\enq
where $N$ $\in$ $L^2(\Gc;\R^d)$ is independent of $(\tilde B,\tilde X_0)$, with zero mean, and unit variance.

\begin{Remark} \label{remIto}
{\rm It\^o's formula \reff{Ito2} is proved in Proposition 6.3 in \cite{cardel14}, and holds true for any function $\tilde u$ which is twice continuously Fr\'echet differentiable on   $L^2(\Gc;\R^d)$. The fact that $\tilde u$ has a lifted structure plays no role, and is used only to derive from \reff{Uu1}-\reff{Uu2} 
It\^o's formula \reff{Ito} on the Wasserstein space $\Pc_{_2}(\R^d)$.  Recall however that It\^o's formula \reff{Ito} holds even if the lift is not twice continuously Fr\'echet differentiable as shown in \cite{chacridel15} (see also \cite{caretal16}). 
\ep
}
\end{Remark}

\subsection{Dynamic programming equation}

The dynamic programming Bellman equation associated to the value function of the stochastic  McKean-Vlasov control problem takes the form:
\begin{equation} \label{HJBdynpro}
\left\{
\begin{array}{rcl}
- \partial_t v - \Inf_{a\in\bA} \Big[ \hat f(\mu,a) + \mu\big(\L^a v(t,\mu) \big) + \mu\otimes\mu\big(\M^a v(t,\mu) \big) \Big] & =& 0,  \;\;\; (t,\mu)  \in [0,T)\times\Pc_{_2}(\R^d), \\
v(T,\mu) &=& \hat g(\mu), \;\;\; \mu \in \Pc_{_2}(\R^d),
\end{array}
\right.
\end{equation}
where for $\phi$ $\in$ $\Cc_b^2(\Pc_{_2}(\R^d))$, $a$ $\in$ $\bA$, and $\mu$ $\in$ $\Pc_{_2}(\R^d)$, $\L^a\phi(\mu)$ $\in$ $L_\mu^2(\R)$ is the function $\R^d$ $\rightarrow$ $\R$ defined by
\beq \label{defL}
\L^a\phi(\mu)(x) &:=& \partial_\mu \phi(\mu)(x).b(x,\mu,a) + \frac{1}{2}{\rm tr}\big(\partial_x\partial_\mu\phi(\mu)(x)( \sigma\sigma\trans + \sigma_0\sigma_0\trans)(x,\mu,a) \big),
\enq
and $\M^a\phi(\mu)$ $\in$ $L_{\mu\otimes\mu}^2(\R)$ is the function $\R^d\times\R^d$ $\rightarrow$ $\R$ defined by 
\beq \label{defM}
\M^a\phi(\mu)(x,x') &:=& \frac{1}{2} {\rm tr}\big(  \partial_\mu^2\phi(\mu)(x,x')\sigma_0(x,\mu,a)\sigma_0\trans(x',\mu,a) \big). 
\enq

Alternatively, by viewing  the value function  as a function on $[0,T]\times L^2(\Gc;\R^d)$ via the lifting identification, and keeping the same notation $v(t,\xi)$ $=$ $v(t,\Lc(\xi))$ (recall that $v$ depends on $\xi$ only via its distribution),  
we see from the connection \reff{Uu1}-\reff{Uu2}  between derivatives in the Wasserstein space $\Pc_{_2}(\R^d)$ and in the Hilbert space $L^2(\Gc;\R^d)$ that 
the Bellman equation \reff{HJBdynpro} is written also in  $[0,T]\times L^2(\Gc;\R^d)$   as 
\begin{equation} \label{HJBdynpro1}
\left\{
\begin{array}{rcl}
- \partial_t v -   H\big (\xi,Dv(t,\xi),D^2 v(t,\xi) \big) & =& 0,  \;\;\; (t,\xi)  \in [0,T)\times L^2(\Gc;\R^d), \\
v(T,\xi) &=& \tilde\E^1\big[ g(\xi,\Lc(\xi)) \big], \;\;\; \xi \in   L^2(\Gc;\R^d),
\end{array}
\right.
\end{equation} 
where $H$ $:$ $L^2(\Gc;\R^d)\times L^2(\Gc;\R^d)\times S(L^2(\Gc;\R^d))$ $\rightarrow$ $\R$ is defined by
\beq 
H(\xi,P,Q) &=&    \Inf_{a\in\bA}   \tilde\E^1\Big[f(\xi,\Lc(\xi),a) + P.b(\xi,\Lc(\xi),a)  \label{defH} \\
& & \;\;\; + \; \frac{1}{2}  Q(\sigma_0(\xi,\Lc(\xi),a)).\sigma_0(\xi,\Lc(\xi),a) +  \frac{1}{2}  Q(\sigma(\xi,\Lc(\xi),a)N).\sigma(\xi,\Lc(\xi),a)N \Big], \nonumber 
\enq
with $N$ $\in$ $L^2(\Gc;\R^n)$ of zero mean, and unit variance, and independent of $\xi$.

The purpose of this section is to prove an analytic characterization of the value function  in terms of the dynamic programming Bellman equation. We shall adopt a notion of viscosity solutions following the approach in \cite{lio12},  
which consists via the lifting identification in working in the Hilbert space $L^2(\Gc;\R^d)$ instead of working in the Wasserstein space $\Pc_{_2}(\R^d)$. Indeed, comparison principles for viscosity solutions in 
the Wasserstein space, or more generally in metric spaces, are difficult to obtain as we have to deal with locally non compact spaces (see e.g. \cite{ambetal05}, \cite{ganetal08}, \cite{fenkat09}), and instead by working in  separable Hilbert spaces, one can essentially reduce to the case of Euclidian spaces by projection, and then take advantage of the results developed  for viscosity solutions, in particular here,  for second order Hamilton-Jacobi-Bellman equations, see 
\cite{lio89b}, \cite{fabgozswi15}.  
We shall assume that the $\sigma$-algebra $\Gc$ is countably generated upto null sets, which ensures that the 
Hilbert space $L^2(\Gc;\R^d)$ is separable, see \cite{doo94}, p. 92.  This is satisfied for example when $\Gc$ is the Borel $\sigma$-algebra of a canonical space $\tilde\Omega^1$ 
of continuous functions on $\R_+$ (see Exercise 4.21 in Chapter 1 of \cite{revyor99}).

\begin{Definition} \label{defvisco}
We say that a continuous function $u$ $:$ $[0,T]\times\Pc_{_2}(\R^d)$ $\rightarrow$ $\R$ is a viscosity (sub, super) solution to \reff{HJBdynpro} if its lifted version $\tilde u$ on $[0,T]\times L^2(\Gc;\R^d)$ is a viscosity 
(sub, super) solution to \reff{HJBdynpro1}, that is: 

\noindent (i) $\tilde u(T,\xi)$ $\leq$ $\tilde\E^1\big[ g(\xi,\Lc(\xi)) \big]$, and for any test function $\varphi$ $\in$ $\Cc^2([0,T]\times L^2(\Gc;\R^d))$ (the set of real-valued continuous functions on $[0,T]\times L^2(\Gc;\R^d)$ which are continuously differentiable in $t$ $\in$ $[0,T)$, and twice continuously Fr\'echet differentiable on $L^2(\Gc;\R^d)$) 
 s.t. $\tilde u-\varphi$ has a maximum at $(t,\xi)$ $\in$ $[0,T)\times L^2(\Gc;\R^d)$, one has
\beqs
- \partial_t \varphi(t,\xi) -  H\big (\xi,D\varphi(t,\xi),D^2 \varphi(t,\xi) \big) &  \leq & 0.
\enqs
\noindent (ii) $\tilde u(T,\xi)$ $\geq$ $\tilde\E^1\big[ g(\xi,\Lc(\xi)) \big]$, and for any test function $\varphi$ $\in$ $\Cc^2([0,T]\times L^2(\Gc;\R^d)$ s.t. $\tilde u-\varphi$ has a minimum at $(t,\xi)$ $\in$ $[0,T)\times L^2(\Gc;\R^d)$, one has
\beqs
- \partial_t \varphi(t,\xi) -  H\big (\xi,D\varphi(t,\xi),D^2 \varphi(t,\xi) \big) &  \geq & 0.
\enqs
\end{Definition} 

\begin{Remark}
{\rm Since the lifted function $\tilde u$ of a smooth solution $u$ $\in$ $\Cc^2([0,T]\times\Pc_{_2}(\R^d))$ to  \reff{HJBdynpro}, may not be smooth in 
$[0,T]\times L^2(\Gc;\R^d)$, it says that $u$ cannot be viewed in general as a viscosity solution to   \reff{HJBdynpro} in the sense of Definition 
\ref{defvisco} unless we add the extra-assumption that its lifted function is indeed twice continuously Fr\'echet differentiable on   $L^2(\Gc;\R^d)$. 
Hence, a more natural and intrinsic definition of viscosity solutions would use test functions on $[0,T]\times\Pc_{_2}(\R^d)$: in this case, it would be possible to get the viscosity property from the dynamic programming principle and It\^o's formula \reff{Ito}, but as pointed out above,  the uniqueness result (and so the characterization) in the Wasserstein space is a challenging issue, beyond the scope of this paper. We have then chosen here to work with test functions on $[0,T]\times L^2(\Gc;\R^d)$, not necessarily of the lifted form.
\ep
}
\end{Remark}

\vspace{1mm}

The main result of this section is the viscosity characterization of the value function for the stochastic McKean-Vlasov control problem  \reff{defv} to the dynamic programming Bellman equation \reff{HJBdynpro} (or \reff{HJBdynpro1}).

\begin{Theorem}
The value function $v$ is the unique continuous viscosity solution to  \reff{HJBdynpro} satisfying a quadratic growth condition \reff{vquadra2}. 
\end{Theorem}
{\bf Proof.}  {\it  (1) Viscosity property}.  Let us first reformulate the dynamic programming principle (DPP)  of Proposition \ref{propdyn}  for the value function viewed now  as a function on $[0,T]\times L^2(\Gc;\R^d)$.  
For this, we take a copy $\tilde B$  of $B$ on  the probability space $(\tilde\Omega^1,\Gc,\tilde\P^1)$, and  given $(t,\xi)$ $\in$ $[0,T]\times L^2(\Gc;\R^d)$, $\alpha$ $\in$ $\Ac$, we consider on 
$(\tilde\Omega=\Omega^0\times\tilde\Omega^1,\tilde\Fc=\Fc^0\otimes\Gc,\tilde\P =\P^0\otimes\tilde\P^1)$ the solution $\tilde X^{t,\xi,\alpha}$, $t\leq s\leq T$, to the McKean-Vlasov equation 
\beqs
\tilde X_s^{t,\xi,\alpha} &=& \xi + \int_t^s b(\tilde X_r^{t,\xi,\alpha},\tilde\P_{\tilde X_s^{t,\xi,\alpha}}^{W^0},\alpha_r) dr  +  
\int_t^s \sigma(\tilde X_r^{t,\xi,\alpha},\tilde\P_{\tilde X_s^{t,\xi,\alpha}}^{W^0},\alpha_r) d\tilde B_r  \\
& & \;\;\;\;\; + \;  \int_t^s \sigma_0(\tilde X_r^{t,\xi,\alpha},\tilde\P_{\tilde X_s^{t,\xi,\alpha}}^{W^0},\alpha_r) dW_r^0, \;\;\; t \leq s \leq T,
\enqs
where $\tilde\P_{\tilde X_s^{t,\xi,\alpha}}^{W^0}$ denotes the regular conditional distribution of $\tilde X_s^{t,\xi,\alpha}$ given $\Fc^0$.  In other words, $\tilde X^{t,\xi,\alpha}$ is a copy of 
$X^{t,\xi,\alpha}$ on $(\tilde\Omega,\tilde\Fc,\tilde\P)$, and denoting by $\check X_s^{t,\xi,\alpha}(\omega^0)$ $=$ 
$\tilde X_s^{t,\xi,\alpha}(\omega^0,.)$, $t\leq s\leq T$, we see that the process $\{\check X_s^{t,\xi,\alpha},t\leq s \leq T\}$ is $\F^0$-progressive, valued in $L^2(\Gc;\R^d)$, and 
$\tilde\P^1_{\check X_s^{t,\xi,\alpha}}$ $=$ $\rho_s^{t,\mu,\alpha}$ for $\mu$ $=$ $\Lc(\xi)$.  Therefore,  the lifted value function on $[0,T]\times L^2(\Gc;\R^d)$ identified with the value function on $[0,T]\times\Pc_{_2}(\R^d)$ satisfies 
$v(s,\check X_s^{t,\xi,\alpha})$ $=$ $v(s,\rho_s^{t,\mu,\alpha})$, $t\leq s \leq T$.  By noting that $\hat f(\rho_s^{t,\mu,\alpha},\alpha_s)$ $=$ $\tilde\E^1\big[ f(\check X_s^{t,\xi,\alpha},\tilde\P_{\check X_s^{t,\xi,\alpha}}^1,\alpha_s)\big]$, 
we  obtain from Proposition \ref{propdyn} the lifted DPP: for all $(t,\xi)$ $\in$ $[0,T]\times L^2(\Gc;\R^d)$, 
\beq
v(t,\xi) &=& \inf_{\alpha\in\Ac} \inf_{\theta\in\Tc_{t,T}^0} 
\E^0\Big[ \int_t^\theta \tilde\E^1\big[ f(\tilde X_s^{t,\xi,\alpha},\tilde\P_{\check X_s^{t,\xi,\alpha}}^1,\alpha_s)\big]  ds \; + \; v(\theta,\check X_\theta^{t,\xi,\alpha}) \Big]
\label{DPPstrong1lift} \\
&=& \inf_{\alpha\in\Ac} \sup_{\theta\in\Tc_{t,T}^0} 
\E^0\Big[ \int_t^\theta \tilde\E^1\big[ f(\tilde X_s^{t,\xi,\alpha},\tilde\P_{\check X_s^{t,\xi,\alpha}}^1,\alpha_s)\big]  ds \; + \; v(\theta,\check X_\theta^{t,\xi,\alpha}) \Big]. \label{DPPstrong2lift}
\enq
We already know that $v$ is continuous on $[0,T]\times L^2(\Gc;\R^d)$, hence in particular at $T$,  so that $v(T,\xi)$ $=$ $\tilde\E^1[g(\xi,\Lc(\xi))]$, and it remains to 
derive the viscosity property for the value function in $[0,T)\times L^2(\Gc;\R^d)$  by following standard arguments that we adapt in our context. 

\vspace{1mm}

\noindent  {\it (i) Subsolution property}.  Fix $(t,\xi)$ $\in$ $[0,T)\times L^2(\Gc;\R^d)$, and consider some test function $\varphi$ $\in$ $\Cc^2([0,T]\times L^2(\Gc;\R^d))$ s.t. $v-\varphi$ has a maximum at $(t,\xi)$, and w.l.o.g. 
$v(t,\xi)$ $=$ $\varphi(t,\xi)$, so that $v$ $\leq$ $\varphi$. Let $a$ be an arbitrary element in  $\bA$, $\alpha$ $\equiv$ $a$ the constant control in $\Ac$ equal to $a$, and consider the stopping time in $\Tc_{t,T}^0$: 
$\theta_{h}$ $=$ $\inf\{s \geq t : \tilde\E^1[|\check X_s^{t,\xi,a} - \xi|^2] \geq \delta^2 \}$ $\wedge$ $(t+h)$,  with $h$ $\in$ $(0,T-t)$, and  $\delta$ some positive constant  small enough (depending on $\xi$), so that 
$\varphi$ and its continuous derivatives $\partial_t\varphi$, $D\varphi$, $D^2\varphi$ are bounded on the ball in $L^2(\Gc;\R^d)$ of center $\xi$ and radius $\delta$. 
From  the first part \reff{DPPstrong1lift} of the DPP, we get
\beqs
\varphi(t,\xi) & \leq & \E^0\Big[ \int_t^{\theta_{h}} \tilde\E^1\big[ f(\tilde X_s^{t,\xi,a},\tilde\P_{\check X_s^{t,\xi,a}}^1,a)\big]  ds \; + \; \varphi(\theta_{h},\check X_{\theta_{h}}^{t,\xi,a}) \Big]. 
\enqs
Applying It\^o's formula \reff{Ito2} to $\varphi(s,\check X_s^{t,\xi,a})$, and noting that the stochastic integral w.r.t. $W^0$ vanishes under expectation $\E^0$ by the localization with the stopping time $\theta_h$, we then have
\beq 
0 & \leq &  \E^0 \Big[ \frac{1}{h} \int_t^{\theta_h} \partial_t\varphi(s,\check X_s^{t,\xi,a}) +  \tilde\E^1 \big[ f(\tilde X_s^{t,\xi,a},\tilde\P_{\check X_s^{t,\xi,a}}^1,a) + D \varphi(s,\check X_s^{t,\xi,a}).b(\check X_s^{t,\xi,a},\tilde\P^1_{\check X_s^{t,\xi,a}},a)   \nonumber \\
& & \;\;\;\;\; + \; \frac{1}{2}  D^2\varphi(s, \check X_s^{t,\xi,a})(\sigma(\check X_s^{t,\xi,a},\tilde\P_{\check X_s^{t,\xi,a}}^1,a)N).\sigma(\check X_s^{t,\xi,a},\tilde\P_{\check X_s^{t,\xi,a}}^1,a)N  \nonumber \\
& &  \;\;\;\;\; + \; \frac{1}{2}  D^2\varphi(s, \check X_s^{t,\xi,a})(\sigma_0(\check X_s^{t,\xi,a},\tilde\P_{\check X_s^{t,\xi,a}}^1,a)).\sigma_0(\check X_s^{t,\xi,a},\tilde\P_{\check X_s^{t,\xi,a}}^1,a) \big] ds \Big] \nonumber \\
& =: &  \E^0 \Big[ \frac{1}{h} \int_t^{\theta_h}  F_s(t,\xi,a) ds  \Big],  \label{Itophi}
\enq
with $N$ $\in$ $L^2(\Gc;\R^n)$ of zero mean, and unit variance, and independent of $(\tilde B,\xi)$.  Since  the map  
 $s$ $\in$ $[t,T]$ $\mapsto$ $\tilde\E^1[\psi(\tilde X_s^{t,\xi,a})]$ $=$ $\E[\psi(X_s^{t,\xi,a})|\Fc^0]$ $=$  $\rho_s^{t,\mu,a}(\psi)$ (for $\mu$ $=$ $\Lc(\xi)$) is continuous $\P^0$-a.s. (recall that $\rho_s^{t,\mu,\alpha}$ is continuous in $s$), 
 for any bounded continuous function $\psi$ on $\R^d$,  we see that the process $\{F_s(t,\xi,a), t\leq s\leq \theta_h\}$ has continuous paths $\P^0$ almost surely.   Moreover, by (standard) It\^o's formula, we have for all $t\leq s \leq T$, 
 \beqs
 \tilde\E^1\big[|\check X_s^{t,\xi,a} - \xi|^2\big]  \; = \;  \E \big[ |X_s^{t,\xi,a} - \xi |^2 | \Fc^0 \big] & = &  \int_t^s \E\big[ 2(X_r^{t,\xi,a} - \xi).b_r + \sigma_r\sigma_r\trans + \sigma_r^0(\sigma_r^0)\trans | \Fc^0  \big] dr \\
 &  &  + \; \int_t^s \E \big[ 2 (X_r^{t,\xi,a} - \xi)\trans\sigma_r^0 | \Fc^0] dW_r^0,  
 \enqs
where we set $b_s$ $=$ $b(X_s^{t,\xi,a},\P_{X_s^{t,\xi,a}}^{W^0},a)$, $\sigma_s$ $=$ $\sigma(X_s^{t,\xi,a},\P_{X_s^{t,\xi,a}}^{W^0},a)$,  
$\sigma_s^0$ $=$ $\sigma_0(X_s^{t,\xi,a},\P_{X_s^{t,\xi,a}}^{W^0},a)$. This shows that the map $s$ $\in$ $[t,T]$ $\mapsto$ $ \tilde\E^1[|\check X_s^{t,\xi,a} - \xi|^2]$  is continuous $\P^0$-a.s., and thus 
$\theta_h(\omega^0)$ $=$ $t+h$ for $h$ small enough ($\leq$ $\bar h(\omega^0)$), $\P^0(d\omega^0)$-a.s.  By the mean-value theorem, we then get $\P^0$ almost surely,  
$\frac{1}{h} \int_t^{\theta_h}  F_s(t,\xi,a) ds$  $\rightarrow$ $F_t(t,\xi,a)$, as $h$ goes to zero, and so from the dominated convergence theorem in \reff{Itophi}:
\beqs
0 \;  \leq \;  F_t(t,\xi,a) & = & \partial_t \varphi(t,\xi) + \tilde\E^1\big[ f(\xi,\Lc(\xi),a)  + D\varphi(t,\xi).b(\xi,\Lc(\xi),a)  \\
& & \;\;\;\;\; + \; \frac{1}{2} D^2\varphi(t,\xi) (\sigma(\xi,\Lc(\xi),a)N).\sigma(\xi,\Lc(\xi),a)N \\
& &  \;\;\;\;\;\;\;  + \;  \frac{1}{2} D^2\varphi(s,\xi) (\sigma_0(\xi,\Lc(\xi),a)).\sigma_0(\xi,\Lc(\xi),a) \Big]. 
\enqs
Since $a$ is arbitrary in $\bA$, this shows the required viscosity subsolution property. 
 
  \vspace{1mm}

\noindent  {\it (ii) Supersolution property}.  Fix $(t,\xi)$ $\in$ $[0,T)\times L^2(\Gc;\R^d)$, and consider some test function $\varphi$ $\in$ $\Cc^2([0,T]\times L^2(\Gc;\R^d))$ s.t. $v-\varphi$ has a minimum at $(t,\xi)$, and w.l.o.g. 
$v(t,\xi)$ $=$ $\varphi(t,\xi)$, so that $v$ $\geq$ $\varphi$.  From the continuity assumptions in {\bf (H1)}-{\bf (H2)}, we observe that the function $\Hc$ defined on $[0,T]\times L^2(\Gc;\R^d)$ by 
\beqs
\Hc(s,\zeta) &:=&  H(\zeta,D\varphi(s,\zeta),D^2\varphi(s,\zeta)), 
\enqs
is continuous.  Then, given an arbitrary  $\eps$ $>$ $0$,  there exists $\bar h$ $\in$ $(0,T-t)$,  $\delta$ $>$ $0$ s.t. for all $s$ $\in$ $[t,t+\bar h]$,  and $\zeta$ $\in$ $L^2(\Gc;\R^d)$ with $\tilde\E^1[|\zeta-\xi|^2]$ $\leq$ $\delta$, 
\beqs
\Big| \big(\partial_t\varphi + \Hc\big)(s,\zeta) -  \big(\partial_t\varphi +  \Hc\big)(t,\xi) \Big| & \leq & \eps. 
\enqs
 From the second part  \reff{DPPstrong2lift} of the DPP, for any $h$ $\in$ $(0,\bar h)$, there exists $\alpha$ $\in$ $\Ac$ s.t. 
\beqs
\varphi(t,\xi) + \eps h  & \geq & \E^0\Big[ \int_t^{\theta_{h}} \tilde\E^1\big[ f(\tilde X_s^{t,\xi,\alpha},\tilde\P_{\check X_s^{t,\xi,\alpha}}^1,\alpha_s)\big]  ds \; + \; \varphi(\theta_{h},\check X_{\theta_{h}}^{t,\xi,\alpha}) \Big], 
\enqs 
where we take $\theta_h$ $=$  $\inf\{s \geq t : \tilde\E^1[|\check X_s^{t,\xi,\alpha} - \xi|^2] \geq \delta^2\}$ $\wedge$ $(t+h)$ (assuming w.l.o.g. that $\delta$ is small enough (depending on $\xi$), so that 
$\varphi$ and its continuous derivatives $\partial_t\varphi$, $D\varphi$, $D^2\varphi$ are bounded on the ball in $L^2(\Gc;\R^d)$ of center $\xi$ and radius $\delta$). 
Applying again It\^o's formula \reff{Ito2} to $\varphi(s,\check X_s^{t,\xi,\alpha})$, and by definition of $\Hc$, we 
get
\beq
\eps & \geq &  \E^0 \Big[ \frac{1}{h} \int_t^{\theta_h} \big(\partial_t\varphi  + \Hc\big)(s,\check X_s^{t,\xi,\alpha}) ds \Big] \nonumber \\
& \geq &   \Big[ \big(\partial_t\varphi +  \Hc\big)(t,\xi) - \eps \Big] \frac{\E^0[\theta_h]-t}{h},  \label{viscosureps} 
\enq
by the choice of $h$, $\delta$, and $\theta_h$. Now, by noting from Chebyshev's inequality that 
\beqs
\P^0[ \theta_h < t+h ] & \leq & \P^0 \big[ \sup_{t\leq s\leq t+h} \tilde\E^1[|\check X_s^{t,\xi,\alpha} - \xi |^2] \geq \delta \big] \\
& \leq & \frac{\E^0 \Big[  \Sup_{t\leq s\leq t+h} \tilde\E^1[|\check X_s^{t,\xi,\alpha} - \xi |^2]  \Big] }{\delta} \; \leq \; \frac{ C(1+ \tilde\E^1[|\xi|^2]) h}{\delta} 
\enqs
and using the   obvious inequality: $1-\P^0[ \theta_h < t+h ]$ $=$ $\P[\theta_h=t+h]$ $\leq$ $\frac{\E^0[\theta_h]-t}{h}$ $\leq$ $1$,  we see that $\frac{\E^0[\theta_h]-t}{h}$ converges to $1$ when $h$ goes to zero, and deduce from 
\reff{viscosureps} that 
 \beqs
 2 \eps & \geq & \big(\partial_t\varphi +  \Hc\big)(t,\xi).  
 \enqs
We obtain the required   viscosity supersolution property by sending $\eps$ to zero.  
 
 \vspace{1mm}
 
\noindent {\it  (2) Uniqueness property}.  In view of our definition of viscosity solution, we have to show a comparison principle for viscosity solutions to the lifted Bellman equation \reff{HJBdynpro1}.    
We use the comparison principle proved in Theorem 3.50 in \cite{fabgozswi15} and only need to check that the hypotheses of this theorem are satisfied in our context for the lifted Hamiltonian $H$ defined in \reff{defH}. 
Notice that the Bellman equation \reff{HJBdynpro1} is a bounded equation in the terminology 
of  \cite{fabgozswi15} (see their section  3.3.1) meaning that there is no linear dissipative operator on $L^2(\Gc;\R^d)$ in the equation. 
Therefore, the notion of $B$-continuity reduces to the standard notion of continuity in $L^2(\Gc;\R^d)$ since one can take for $B$ the identity operator. 
Their Hypothesis 3.44 follows from the uniform continuity  of $b$, $\sigma$, $\sigma_0$ and $f$ in {\bf (H1)}-{\bf (H2)}.  Hypothesis 3.45  is immediately satisfied since there is no discount factor in our equation, i.e. 
$H$ does not depend on $v$ but only on its derivatives.  
The monotonicity condition in $Q$ $\in$ $S(L^2(\Gc;\R^d))$ of $H$ in Hypothesis 3.46  is clearly satisfied.  Hypothesis 3.47  holds directly  when dealing with  bounded equations.   
Hypothesis 3.48  is obtained from the Lipschitz condition of $b,\sigma,\sigma_0$ in {\bf (H1)}, and the uniform continuity condition on $f$ in {\bf (H2)}, while Hypothesis 3.49  follows from the growth condition of 
$\sigma$, $\sigma_0$ in  {\bf (H1)}. One can then apply Theorem 3.50  in  \cite{fabgozswi15} and conclude that comparison principle holds for the Bellman equation \reff{HJBdynpro1}.  
\ep

\vspace{2mm}

We conclude this section with a verification theorem, which gives an analytic feedback form  of the optimal control when there is a smooth solution to the Bellman equation \reff{HJBdynpro} in the Wasserstein space. 
We refer to the recent paper \cite{ganswi16} for existence result of smooth solution to the Bellman equation on small time horizon.

\begin{Theorem} \label{theoverif} (Verification theorem)

\noindent Let $w$ $:$ $[0,T]\times\Pc_{_2}(\R^d)$ $\rightarrow$ $\R$ be a function in $\Cc_b^{1,2}([0,T]\times\Pc_{_2}(\R^d))$, i.e.   $w$ is continuous on $[0,T]\times\Pc_{_2}(\R^d)$,  $w(t,.)$ $\in$ $\Cc_b^2(\Pc_{_2}(\R^d))$, and 
$w(.,\mu)$ $\in$ $C^1([0,T))$, and satisfying a quadratic growth condition as in  \reff{vquadra2}, together with a linear growth condition for its derivative:
\beq \label{Dwmulin}
|\partial_\mu w(t,\mu)(x)| & \leq & C(1 + |x| + \|\mu\|_{_2}), \;\;\; \forall (t,x,\mu) \in [0,T]\times\R^d\times\Pc_{_2}(\R^d), 
\enq
for some positive constant $C$. 
Suppose that  $w$ is solution to the Bellman equation \reff{HJBdynpro}, and there exists  for all $(t,\mu)$ $\in$ $[0,T)\times\Pc_{_2}(\R^d)$ an element $\hat a(t,\mu)$ $\in$ $\bA$ attaining the infimum in \reff{HJBdynpro} s.t.
the map $(t,\mu)$ $\mapsto$ $\hat a(t,\mu)$ is measurable, and the stochastic McKean-Vlasov equation
\beqs
d\hat X_s &=& b(\hat X_s,\P_{\hat X_s}^{W^0},\hat a(s,\P_{\hat X_s}^{W^0})) ds +   \sigma(\hat X_s,\P_{\hat X_s}^{W^0},\hat a(s,\P_{\hat X_s}^{W^0})) dB_s \\
& & \;\;\;+ \;  \sigma(\hat X_s,\P_{\hat X_s}^{W^0},\hat a(s,\P_{\hat X_s}^{W^0})) dW_s^0, \;\;\; t \leq s \leq T,  \; \hat X_t \; = \xi, 
\enqs
admits a unique solution denoted $(\hat X_s^{t,\xi})_{t\leq s \leq T}$, for any $(t,\xi)$ $\in$ $[0,T]\times L^2(\Gc;\R^d)$ (This is satisfied e.g. when $\mu$ $\mapsto$ $\hat a(t,\mu)$ is Lipschitz on $\Pc_{_2}(\R^d)$).  
Then,  $w$ $=$ $v$, and the feedback control $\alpha^*$ $\in$ $\Ac$ defined by
\beq \label{defalphaopt}
\alpha_s^* &=& \hat a(s,\P^{W^0}_{\hat X_s^{t,\xi}}), \;\;\; t \leq s <  T,
\enq
is an optimal control for $v(t,\mu)$, i.e.  $v(t,\mu)$ $=$ $J(t,\mu,\alpha^*)$, with $\mu$ $=$ $\Lc(\xi)$. 
\end{Theorem}
{\bf Proof.}   Fix $(t,\mu=\Lc(\xi))$ $\in$ $[0,T]\times\Pc_{_2}(\R^d)$, and consider some arbitrary control $\alpha$ $\in$ $\Ac$ associated to $\rho_s^{t,\mu,\alpha}$ $=$ $\P_{X_s^{t,\xi,\alpha}}^{W^0}$, $t\leq s\leq T$. 
Denote by $X_s^{'t,\xi,\alpha}$ a copy of $X_s^{t,\xi,\alpha}$ on another probability space 
$(\Omega'=\Omega^0\times\Omega^{'1}\,\Fc^0\otimes\Fc^{'1},\P^0\times\P^{'1})$, with  $(\Omega^{'1},\Fc^{'1},\P^{'1})$ supporting  $B'$ a copy of $B$.   Applying It\^o's formula \reff{Ito} to 
$w(s,\rho_s^{t,\mu,\alpha})$ between $t$ and the $\F^0$-stopping time $\theta_T^n$ $=$ $\inf\{s\geq t: \|\rho_s^{t,\mu,\alpha}\|_{_2} \geq n\} \wedge T$, we obtain
\beq
& & w(\theta_T^n,\rho_{\theta_T^n}^{t,\mu,\alpha}) \nonumber \\
&=& w(t,\mu) +  \int_t^{\theta_T^n} \Big\{ \Dt{w}(s,\rho_s^{t,\mu,\alpha})  
+  \E_{_{W^0}} \Big[ \partial_\mu w(s,\rho_s^{t,\mu,\alpha})(X_s^{t,\xi,\alpha}).b(X_s^{t,\xi,\alpha},\rho_s^{t,\mu,\alpha},\alpha_s) \nonumber \\
& & + \; \frac{1}{2} {\rm tr}\big[ \partial_x\partial_\mu w(s,\rho_s^{t,\mu,\alpha})(X_s^{t,\xi,\alpha})(\sigma\sigma\trans(X_s^{t,\xi,\alpha},\rho_s^{t,\mu,\alpha},\alpha_s) +  
\sigma_0\sigma_0\trans(X_s^{t,\xi,\alpha},\rho_s^{t,\mu,\alpha},\alpha_s) ) \big]  \Big] \nonumber \\
& & + \; \E_{_{W^0}} \Big[ \E_{_{W^0}}' \big[  \frac{1}{2} {\rm tr}\big( \partial_\mu^2 w(s,\rho_s^{t,\mu,\alpha})(X_s^{t,\xi,\alpha},X_s^{'t,\xi,\alpha}) \sigma_0(X_s^{t,\xi,\alpha},\rho_s^{t,\mu,\alpha},\alpha_s)
\sigma_0\trans(X_s^{'t,\xi,\alpha},\rho_s^{t,\mu,\alpha},\alpha_s) \big) \big]  \Big] \Big\} ds \nonumber \\
& & + \; \int_t^{\theta_T^n} \E_{_{W^0}} \big[ \partial_\mu w(s,\rho_s^{t,\mu,\alpha})(X_s^{t,\mu,\alpha})\trans\sigma_0(X_s^{t,\mu,\alpha},\rho_s^{t,\mu,\alpha},\alpha_s) \big] dW_s^0 \nonumber \\
&=& w(t,\mu) +  \int_t^{\theta_T^n} \Big[ \Dt{w}(s,\rho_s^{t,\mu,\alpha})  + \rho_s^{t,\mu,\alpha} \big(  \L^{\alpha_s} w(s,\rho_s^{t,\mu,\alpha}) \big) 
+  \rho_s^{t,\mu,\alpha}\otimes\rho_s^{t,\mu,\alpha}\big( \M^{\alpha_s} w(s,\rho_s^{t,\mu,\alpha}) \big) \Big] ds \nonumber \\
& & + \; \int_t^{\theta_T^n}  \E_{_{W^0}} \big[ \partial_\mu w(s,\rho_s^{t,\mu,\alpha})(X_s^{t,\mu,\alpha})\trans\sigma_0(X_s^{t,\mu,\alpha},\rho_s^{t,\mu,\alpha},\alpha_s) \big] dW_s^0, \label{Itoverif}
\enq
by definition of $\L^a$ and $\M^a$ in \reff{defL}-\reff{defM}, and recalling again that $\rho_s^{t,\mu,\alpha}$ $=$ $\P_{X_s^{t,\xi,\alpha}}^{W^0}$.  Now, the integrand  of the stochastic integral w.r.t. $W^0$ in \reff{Itoverif} satisfies:
\beqs
& & \Big| \E_{_{W^0}} \big[ \partial_\mu w(s,\rho_s^{t,\mu,\alpha})(X_s^{t,\mu,\alpha})\trans\sigma_0(X_s^{t,\mu,\alpha},\rho_s^{t,\mu,\alpha},\alpha_s) \big] \Big|^2 \\
& \leq &   \Big(  \int_{\R^d} \big| \partial_\mu w(s,\rho_s^{t,\mu,\alpha})(x)\trans\sigma_0(x,\rho_s^{t,\mu,\alpha},\alpha_s) \big| \rho_s^{t,\mu,\alpha}(dx) \Big)^2 \\
& \leq & \int_{\R^d} \big| \partial_\mu w(s,\rho_s^{t,\mu,\alpha})(x)\big|^2 \rho_s^{t,\mu,\alpha}(dx)   \int_{\R^d} \big| \sigma_0(x,\rho_s^{t,\mu,\alpha},\alpha_s)  \big|^2 \rho_s^{t,\mu,\alpha}(dx)  \\
& \leq & C(1 + n^2)^2 \; < \; \infty,  \;\;\;\;\;  t \leq s \leq \theta_T^n,
%\Sup_{\|\pi\|_{_2}\leq n} \int_{\R^d} \big| \partial_\mu w(s,\pi)(x)\big|^2 \pi(dx)  \; <  \infty, \;\;\;\;\;  t \leq s \leq \theta_T^n, 
\enqs
from Cauchy-Schwarz inequality, the linear growth condition of $\sigma_0$ in {\bf (H1)}, the choice of  $\theta_T^n$, and condition \reff{Dwmulin}. 
Therefore, the stochastic integral in \reff{Itoverif} vanishes in $\E^0$-expectation, and  we get
\beq 
\E^0 \big[ w(\theta_T^n,\rho_{\theta_T^n}^{t,\mu,\alpha}) \big] &=& w(t,\mu) \; + \; \E^0 \Big[ \int_t^{\theta_T^n}  \Dt{w}(s,\rho_s^{t,\mu,\alpha})  + \rho_s^{t,\mu,\alpha} \big(  \L^{\alpha_s} w(s,\rho_s^{t,\mu,\alpha}) \big)  \nonumber  \\
& & \hspace{3cm}  + \;   \rho_s^{t,\mu,\alpha}\otimes\rho_s^{t,\mu,\alpha}\big( \M^{\alpha_s} w(s,\rho_s^{t,\mu,\alpha}) \big)  ds \Big]  \nonumber  \\
& \geq &  w(t,\mu) \; -  \; \E^0 \Big[ \int_t^{\theta_T^n}  \hat f(\rho_s^{t,\mu,\alpha},\alpha_s) ds \Big],  \label{Itoverif2}
\enq
since $w$ satisfies the Bellman equation \reff{HJBdynpro}.  By sending $n$ to infinity into \reff{Itoverif2}, and from the dominated convergence theorem (under the condition that $w$, $f$ satisfy a quadratic growth condition and recalling the estimation \reff{estimrho0}), we obtain:
\beqs
w(t,\mu) & \leq & J(t,\mu,\alpha) \; = \; \E^0\Big[ \int_t^T \hat f(\rho_s^{t,\mu,\alpha},\alpha_s) ds  + \hat g(\rho_T^{t,\mu,\alpha})  \Big].
\enqs
Since $\alpha$ is arbitrary in $\Ac$, this shows that $w$ $\leq$ $v$.

Finally, by applying the same It\^o's argument  with the feedback control $\alpha^*$ $\in$ $\Ac$ in \reff{defalphaopt}, and noting that $\hat X_s^{t,\xi}$ $=$ $X_s^{t,\xi,\alpha^*}$, $\P_{\hat X_s^{t,\xi}}^{W^0}$ $=$ $\rho_s^{t,\mu,\alpha^*}$, 
we have now equality in \reff{Itoverif2}, hence $w(t,\mu)$ $=$ $J(t,\mu,\alpha^*)$ ($\geq$ $v(t,\mu)$), and thus finally the required equality:  $w(t,\mu)$ $=$ $v(t,\mu)$ $=$ $J(t,\mu,\alpha^*)$. 
\ep

\section{Linear quadratic stochastic McKean-Vlasov control}

\setcounter{equation}{0} \setcounter{Assumption}{0}
\setcounter{Theorem}{0} \setcounter{Proposition}{0}
\setcounter{Corollary}{0} \setcounter{Lemma}{0}
\setcounter{Definition}{0} \setcounter{Remark}{0}

We consider the linear-quadratic (LQ) stochastic McKean-Vlasov control problem where the control set $\bA$ is a functional space, which corresponds to  the McKean-Vlasov problem with common noise as presented in the introduction. 
%The case of an Euclidian control set, corresponding to partial observation control problem,  will be addressed in the paper \cite{phawei16}.  

%\subsection{LQ McKean-Vlasov with functional control set} \label{LQfonction}

The control set $\bA$ is the set $L(\R^d;\R^m)$   of Lipschitz functions from $\R^d$ into $A$ $=$ $\R^m$, and  
we consider a multivariate linear McKean-Vlasov controlled dynamics with coefficients given by
\begin{equation} \label{bsigLQ}
%\left\{
\begin{array}{ccl}
b(x,\mu,a) &=&  b_0 + B x + \bar B \bar\mu + C  a(x), \\
\sigma(x,\mu,a) &=& \vartheta +  D x + \bar D \bar\mu + F  a(x),\\
\sigma_0(t, x, \mu, a) &=& \vartheta_0 + D_0 x + \bar D_0 \bar\mu + F_0  a(x),
\end{array}
%\right.
\end{equation}
for $(x,\mu,a)$ $\in$ $\R^d\times\Pc_{_2}(\R^d)\times L(\R^d;\R^m)$, where we set
\beqs
\bar\mu & := & \int_{\R^d} x \mu(dx).  
\enqs
%\beqs
%dX_t &=& \big(B(t)X_s + \bar B(t) \E[X_t] + C(t) \alpha_t + \bar C(t) \E[\alpha_t] \big) dt  \\
%& & \;\;\; + \; \big(  D(t) X_t + \bar D(t) \E[X_t] + F(t) \alpha_t  + \bar F(t) \E[\alpha_t] \big) dB_t,
%\enqs
Here  $B$, $\bar B$,  $D$, $\bar D$, $D_0$, $\bar D_0$,  are constant matrices in $\R^{d\times d}$, $C$,  $F$, $F_0$ are constant matrices in $\R^{d\times m}$, and $b_0$, $\vartheta$, $\vartheta_0$ are constant vectors  in $\R^d$.
The quadratic cost functions  are given by 
\begin{equation} \label{fgLQ}
%\left\{
\begin{array}{rcl}
f(x,\mu,a) &=& x\trans Q_2 x + \bar\mu\trans \bar Q_2 \bar \mu  + a(x)\trans R_2 a(x)  
%+ 2 x\trans M_2 a(x) 
\\
%& &  \;   + \;   q_1.x + \bar q_1.\bar\mu + r_1.a(x), \\
g(x,\mu) &=&  x\trans P_2 x + \bar\mu\trans \bar P_2 \bar\mu, 
%+ p_1.x + \bar p_1.\bar\mu,
\end{array}
%\right.
\end{equation}
where $Q_2$, $\bar Q_2$, $P_2$, $\bar P_2$ are constant matrices  in $\R^{d\times d}$, $R_2$ is a constant matrix in $\R^{m\times m}$. 
%$M_2$ is a constant matrix  in $\R^{d\times m}$. 
%$q_1$, $\bar q_1$,  $p_1$, $\bar p_1$ are constants vectors  in $\R^d$, and $r_1$, $\bar r_1$ are constant vectors in $\R^m$.  
Since $f$ and $g$ are real-valued, we may assume w.l.o.g. that all the matrices $Q_2$, $\bar Q_2$,  $R_2$, $P_2$, $\bar P_2$ are symmetric. We denote by $\S^d$ the set of symmetric matrices in $\R^{d\times d}$, by $\S_+^d$ the subset of nonnegative symmetric matrices, by $\S_{>+}^d$ the subset of symmetric  positive definite matrices, and similarly for $\S^m$, $\S^m_+$, $\S_{>+}^m$.

The functions $\hat f$ and $\hat g$ defined in \reff{hatfg} are then given by
\begin{equation} \label{hatfgLQ1}
\left\{
\begin{array}{rcl}
\hat f(t,\mu,a) &=& {\rm Var}(\mu)(Q_2) + \bar\mu \trans(Q_2 + \bar Q_2)\bar\mu +  \overline{a\star\mu_{_2}}(R_2)  \\
% & &  \;\;+\;\;{\rm Var}(a \star \mu)(R_2) + \overline{a\star\mu} \trans (R_2) \overline{a \star \mu} \; + \;2 \bar \mu \trans M_2\overline{a \star \mu} \\
 %\; + \;  2 \int_{\R^d}  (x-\bar \mu)\trans M_2 a(x) \mu(dx) \\
% & & \;\; +\;\;\big(q_1 + \bar q_1\big).\bar\mu +  r_1.\overline{a \star\mu} \\
 \hat g(\mu) &=& {\rm Var}(\mu)(P_2) + \bar\mu \trans(P_2 + \bar P_2)\bar\mu 
%+ (p_1 + \bar p_1).\bar\mu,
\end{array}
\right.
\end{equation}
for any $\mu$ $\in$ $\Pc_{_2}(\R^d)$, $a$ $\in$ $\bA$ $=$ $L(\R^d;\R^m)$, 
where  we set for any  $\Lambda$ in $\S^{d}$ (resp. in $\S^{m}$), and
$\mu$ $\in$ $\Pc_{_2}(\R^d)$ (resp. $\Pc_{_2}(\R^m)$):
\beqs
\bar\mu_{_2}(\Lambda)  \; := \;  \int x\trans \Lambda x \mu(dx),  & & {\rm Var}(\mu)(\Lambda) \; := \; \bar\mu_{_2}(\Lambda)  - \bar\mu\trans\Lambda\bar\mu,
\enqs
and $a\star\mu$ $\in$ $\Pc_{_2}(\R^m)$ is the image by $a$ $\in$ $L(\R^d;\R^m)$ of the measure $\mu$ $\in$ $\R^m$, so that 
\beqs
\overline{a \star \mu} \; = \; \int_{\R^d} a(x)\mu(dx),  & & \overline{a \star \mu}_{_2}(\Lambda)  \; := \;  \int a(x) \trans \Lambda a(x) \mu(dx). 
\enqs

We look for a value function solution to the Bellman equation \reff{HJBdynpro} in the form
\beq \label{wquadra1}
w(t,\mu) &=& {\rm Var}(\mu)(\Lambda(t)) + \bar\mu\trans\Gamma(t)\bar\mu + \bar\mu\trans\gamma(t) + \chi(t),
\enq
for some  functions $\Lambda$, $\Gamma$ $\in$ $C^1([0,T];\S^d)$, $\gamma$ $\in$ $C^1([0,T];\R^d)$, and $\chi$ $\in$ $C^1([0,T];\R)$.  One easily checks   that 
$w$ lies in $\Cc_b^{1,2}([0,T]\times\Pc_{_2}(\R^d))$ with
%The lift function of $w$ in \reff{wquadra1} is given by
%\beqs
%\tilde w(t,\xi) &=& \tilde\E^1[\xi\trans\Lambda(t)\xi] + \tilde\E^1[\xi]\trans(\Gamma(t)-\Lambda(t))\tilde\E^1[\xi]  + \gamma(t).\tilde\E^1[\xi] + \chi(t),
%\enqs
%for $\xi$ $\in$ $L^2(\Gc;\R^d)$, from which we easily see that $\tilde w$ is Fr\'echet differentiable with Fr\'echet derivative in $L^2(\Gc;\R^d)$ given by: 
%$D\tilde w(t,\xi)$ $=$ $2\xi\trans \Lambda(t) +  2 \tilde\E^1[\xi]\trans(\Gamma(t) - \Lambda(t)) + \gamma(t)$. 
%By computing for all $Y$ $\in$ $L^2(\Fc;\R^d)$ the difference
%\beqs
%\Wc(t, X+Y)-\Wc(t, X) &=&
%\E[(X+Y)\trans \Lambda(t) (X+Y)] -\E[X \trans \Lambda(t) X]\\
% & & \;\;- \; \E[X+Y]\trans(\Lambda(t)  - \Gamma(t))\E[X+Y] + \E[X]\trans(\Lambda(t) - \Gamma(t))\E[X]\\
% & &\;\; +\; \gamma(t). (\E[X+Y]-\E[X])\\
%&=& \E[Y \trans \Lambda(t) X + X \trans \Lambda(t) Y] -  \E[Y] \trans (\Lambda(t) - \Gamma(t)) \E[X] -  \E[X]\trans (\Lambda(t) - \Gamma(t))\E[Y] \\
%& &\;\; +\; \gamma(t) . \E[Y] + o(\vert Y \vert_{L^2}) \\
%&=&2 \E[X \trans \Lambda(t) Y] (+?) \red{-}   2\E[X] \trans(\Lambda(t)- \Gamma(t)) \E[Y] + \gamma(t) \E[Y] + o(\vert Y \vert_{L^2})\\
%&=&
%\E\Big[\big(2X \trans \Lambda(t) +  2 \E[X]\trans(\Gamma(t) - \Lambda(t)) + \gamma(t)\big)Y\Big] \; + \; o(  \|Y \|_{_{L^2}}),
%\enqs
%we see that $\Wc$ is Fr\'echet differentiable (w.r.t. $X$) with $[D\Wc](t,X)(Y)$ $=$ $\E\big[\big(2X \trans \Lambda(t) +  2 \E[X]\trans(\Gamma(t) - \Lambda(t)) + \gamma(t)\big)Y\big]$.
%This shows that $w$ lies in $\Cc_b^{1,2}([0,T]\times\Pc_{_2}(\R^d))$ with
\beqs
\partial_t w(t, \mu)  &=& {\rm Var}(\mu) (\Lambda'(t)) +\bar\mu\trans\Gamma'(t)\bar\mu + \gamma'(t)\bar\mu + \chi'(t), \\
\partial_\mu w(t, \mu)(x)&=& 2 \Lambda(t)(x-\bar\mu) + 2 \Gamma(t)\bar\mu + \gamma(t), \\
\partial_x\partial_\mu w(t, \mu)(x)&=& 2 \Lambda(t),\\
\partial_\mu^2 w(t, \mu)(x,x') &=& 2(\Gamma(t)-\Lambda(t)).
\enqs
Together with the quadratic expression \reff{hatfgLQ1} of $\hat f$, $\hat g$, we then see after some tedious but direct calculations that $w$ satisfies the Bellman equation \reff{HJBdynpro} iff
\beq
& & {\rm Var}(\mu)(\Lambda(T)) + \bar\mu\trans\Gamma(T)\bar\mu + \bar\mu\trans\gamma(T) + \chi(T) \nonumber \\
&=&  {\rm Var}(\mu)(P_2) + \bar\mu \trans(P_2 + \bar P_2)\bar\mu, \label{HJBLQT1} 
%+ (p_1 + \bar p_1).\bar\mu, 
\enq
holds for all  $\mu$ $\in$ $\Pc_{_2}(\R^d)$, and
\beq
& & {\rm Var} (\mu)\big(\Lambda'(t)+ Q_2 + D\trans \Lambda(t)D  + D_0\trans \Lambda(t) D_0  + \Lambda(t) B + B\trans\Lambda(t) \big) +  \inf_{a \in L(\R^d; \R^m)} G_t^{\mu}(a) \nonumber \\
& & \;\;+\; \bar\mu\trans\Big(\Gamma'(t) +Q_2 +\bar Q_2 + (D + \bar D)\trans \Lambda(t)(D +\bar D)\nonumber\\
& & \;\;\;\;\;+ \; (D_0 + \bar D_0) \trans \Gamma(t)(D_0 + \bar D_0)   + \;  \Gamma(t)(B+ \bar B) + (B \;+\; \bar B)\trans \Gamma(t) \Big) \bar\mu \nonumber \\
& & \;  +  \;  \bar\mu\trans\big(\gamma'(t)  + (B+\bar B)\trans\gamma(t)    + 2 (D + \bar D)\trans \Lambda(t) \vartheta  + 2  (D_0 + \bar D_0)\trans\Gamma(t) {\vartheta_0}    + 2 \Gamma(t) b_0  \big)  \nonumber\\
& &  \;\;\;+ \;  \chi'(t)   \; + \;  \gamma(t)\trans b_0 + \vartheta\trans \Lambda(t)\vartheta + {\vartheta_0} \trans \Gamma(t) \vartheta_0 \nonumber \\ 
& =&  0, \label{HJBLQt1}
%& &+(\gamma'(t)+ q_1(t) +\bar q_1(t)+ \gamma(t)(B(t)+ \bar B(t)))\bar\mu +\chi'(t) =0
\enq
holds for all  $t$ $\in$ $[0,T)$, $\mu$ $\in$ $\Pc_{_2}(\R^d)$, where the function $G_t^\mu$ $:$ $L(\R^d;\R^m)$ $\rightarrow$ $\R$ is defined by
\beqs
G_t^\mu(a)&=& {\rm Var}(a \star \mu)(U_t) \;+\; \overline{a \star \mu}\trans V_t   \overline{a \star \mu}
\; +\; 2  \int_{\R^d} (x - \bar \mu)\trans S_t  a(x)\mu(dx)  \nonumber \\
 & & \;+\; 2\bar\mu\trans Z_t \overline{a \star \mu}   \;+\;  Y_t . \overline{a \star \mu}, \label{defGk}
 \enqs
 and we set  $U_t$ $=$ $U(t,\Lambda(t))$, $V_t$ $=$ $V(t,\Lambda(t),\Gamma(t))$, $S_t$ $=$ $S(t,\Lambda(t))$, $Z_t$ $=$ $Z(t,\Lambda(t),\Gamma(t))$,
 $Y_t$ $=$ $Y(t,\Gamma(t),\gamma(t))$ with
\begin{equation} \label{defUVSZ}
\left\{
\begin{array}{rcl}
 U(t,\Lambda(t)) &=&F\trans \Lambda(t) F+ F_0 \trans \Lambda(t)  F_0 +  R_2,\\
V(t,\Lambda(t),\Gamma(t)) &=& F\trans \Lambda(t) F + F_0\trans \Gamma(t) F_0  + R_2\\
%& &\;\;\; + R_2(t) + \bar R_2(t),\\
 S(t,\Lambda(t)) &=&D\trans \Lambda(t) F + D_0 \trans \Lambda(t) F_0 + \Lambda(t) C + M_2,\\
 Z(t,\Lambda(t),\Gamma(t)) &=& (D + \bar D) \trans \Lambda(t)F + (D_0+\bar D_0) \trans\Gamma(t)F +  \Gamma(t)C + M_2 \\
 Y(t,\Gamma(t),\gamma(t)) &=&  C \trans \gamma(t)   + 2 F \trans  \Lambda(t) \vartheta + 2 F_0 \trans  \Gamma(t) \vartheta_0.
\end{array}
\right.
\end{equation}
%Here, $L^2(\mu;\R^m)$ $\supset$ $L(\R^d;\R^m)$ is the Hilbert space of measurable functions on  $\R^d$ valued in $\R^m$ and square integrable w.r.t. $\mu$ $\in$ $\Pc_{_2}(\R^d)$.
 %We now search for the infimum of the function $G_t^{\mu}$, and shall make the following assumptions on the  symmetric matrices of the quadratic cost functional:
 %\begin{equation} \label{condmatrice}
%\left\{
%\begin{array}{ccl}
 %P_2  \; \geq \; 0, \; P_2 +\bar P_2  \; \geq \; 0, & & Q_2(t) \; \geq \; 0, \; Q_2(t) + \bar Q_2(t) \; \geq \; 0,  \\
 %& & R_2(t) \; >  \; 0, \; R_2(t) + \bar R_2(t) \; >  \;0
%\end{array}
%\right.
%\end{equation}
Then, under the condition that the symmetric matrices $U_t$ and $V_t$ in \reff{defUVSZ} are positive, hence invertible (this will be discussed later on), we get after square completion:
\beqs
G_t^\mu(a) &=& {\rm Var}((a-a^*(t,.,\mu))\star \mu)(U_t) \;+\; \overline{ (a-a^*(t,.,\mu)) \star \mu}\trans V_t   \overline{(a-a^*(t,.,\mu)) \star \mu} \\
& & \; - \;   {\rm Var}(\mu)\big(  S_t  U_t^{-1} S_t \trans \big)    -   \bar\mu\trans \big( Z_t V_t^{-1}Z_t\trans\big) \bar\mu    -   Y_t\trans V_t^{-1}Z_t \trans \bar\mu   -  \frac{1}{4}  Y_t\trans V_t^{-1} Y_t.
\enqs
where $a(t,.,\mu)$ $\in$ $L^(\R^d;\R^m)$ is given by
\beq
a^*(t,x,\mu)  &=&   - U_t^{-1} S_t\trans(x-\bar\mu)\;  -\;   V_t^{-1} Z_t \trans \bar\mu \;-\;
\frac{1}{2} V_t^{-1}Y_t.   \label{optialphaLQ}
\enq
This means that $G_t^\mu$ attains its infimum  at   $a^*(t,.,\mu)$, and plugging the above expression of  $G_t^\mu(a^*(t,.,\mu))$ in \reff{HJBLQt1}, we observe that  the relation  \reff{HJBLQT1}-\reff{HJBLQt1}, hence the Bellman equation,
is satisfied by identifying the terms in ${\rm Var}(.)$, $\bar\mu\trans(.)\bar\mu$,  $\bar\mu$,  which leads to the system of ordinary differential
equations (ODEs) for  $(\Lambda,\Gamma,\gamma,\chi)$:
\begin{equation} \label{Riccatilambda1}
\left\{
\begin{array}{rcl}
\Lambda'(t)+ Q_2 + D\trans \Lambda(t)D +D_0 \trans \Lambda(t) D_0 +\Lambda(t)B + B \trans \Lambda(t) && \\
 \;\;\; -S(t,\Lambda(t))U(t,\Lambda(t))^{-1} S(t,\Lambda(t))\trans &= & 0,\\
\Lambda(T) &=& P_2,
\end{array}
\right.
\end{equation}
\begin{equation}\label{Riccatigamma1}
\left\{
\begin{array}{rcl}
\Gamma'(t) +Q_2 +\bar Q_2 + (D + \bar D) \trans \Lambda(t)(D +\bar D) & &  \\
\;\;\; +(D_0 +\bar D_0 ) \trans \Gamma(t)(D_0 + \bar D_0)+ \Gamma(t) \trans(B + \bar B )&& \\
\;\;+\; (B + \bar B)\trans \Gamma(t)- Z(t,\Lambda(t),\Gamma(t))  V(t,\Lambda(t),\Gamma(t))^{-1} Z(t,\Lambda(t),\Gamma(t))\trans &=& 0,\\
\Gamma(T) \; &=& \;  P_2 + \bar P_2,
\end{array}
\right.
\end{equation}
\begin{equation} \label{lingamma1}
\left\{
\begin{array}{rcl}
\gamma'(t) +  \big(B + \bar B)\trans \gamma(t)  -  Z(t,\Lambda(t),\Gamma(t)) V(t,\Lambda(t),\Gamma(t))^{-1} Y(t,\Gamma(t),\gamma(t))   \\
\;\;\; + \; 2 \big( D + \bar D \big)\trans \Lambda(t) \vartheta +2\big(D_0+\bar D_0 \big) \trans \Gamma(t) \vartheta_0 + 2 \Gamma(t) b_0  &=& 0,\\
\gamma(T) \;& =& \;  0
\end{array}
\right.
\end{equation}
\begin{equation} \label{linchi1}
\left \{
\begin{array}{rcl}
\chi'(t) -\; \frac{1}{4} Y(t,\Gamma(t),\gamma(t))\trans V(t,\Lambda(t),\Gamma(t))^{-1}Y(t,\Gamma(t),\gamma(t)) \\
 \;\;\; +  \; \gamma(t)\trans b_0  + \vartheta\trans\Lambda(t)\vartheta + \vartheta_0 \trans \Gamma(t)\vartheta_0 &=& 0, \\
\chi(T) \; &=& \;  0.
\end{array}
\right.
\end{equation}
Therefore, the resolution of the Bellman equation in the LQ framework is reduced to the resolution of the Riccati equations  \reff{Riccatilambda1}  and  \reff{Riccatigamma1} for $\Lambda$ and $\Gamma$, and then
given $(\Lambda,\Gamma)$, to the resolution  of the  linear ODEs \reff{lingamma1} and \reff{linchi1} for $\gamma$ and $\chi$.  Suppose that there exists a solution $(\Lambda,\Gamma)$ $\in$
$C^1([0,T];\S^d)\times C^1([0,T];\S^d)$  to \reff{Riccatilambda1}-\reff{Riccatigamma1}
s.t. $(U_t,V_t)$ in \reff{defUVSZ} lies in $\S^m_{>+}\times\S^m_{>+}$  for all $t$ $\in$ $[0,T]$ (see Remark \ref{remriccati}). Then,  the above calculations
are justified a posteriori, and by noting also that  the mapping $(x,\mu)$ $\mapsto$ $a^*(t,x,\mu)$ is Lipschitz on   $\R^d\times\Pc_{_2}(\R^d)$, 
we deduce by the verification theorem that the value function $v$ is equal to $w$ in \reff{wquadra1} with $(\Lambda,\Gamma,\gamma,\chi)$ solution to
\reff{Riccatilambda1}-\reff{Riccatigamma1}-\reff{lingamma1}-\reff{linchi1}.  Moreover, the optimal control is given in feedback form from \reff{optialphaLQ} by
\beq 
\alpha_t^*(X_t^*)  & = &  a^*(t,X_t^*,\P_{X_t^*}^{W^0}) \nonumber \\
 & =&  - U_t^{-1} S_t\trans\big(X_t^* - \E[X_t^*|\Fc_t^0] \big) -   V_t^{-1} Z_t \trans \E[X_t^*|\Fc_t^0] \;-\;  \frac{1}{2} V_t^{-1}Y_t, \label{opticontrolLQ}
\enq
where $X^*$ is the state process controlled by $\alpha^*$.

\begin{Remark} \label{remriccati}
{\rm 
%In the case where $M_2$ $=$  $0$  (i.e. no crossing term between the state and the control in the quadratic cost function $f$), it is shown in  Proposition 3.1 and 3.2 in  
It is known from \cite{won68}  that  under the condition 
\beq \label{condmatrice}
%\left\{
 P_2  \; \geq \; 0, \; P_2 +\bar P_2  \; \geq \; 0, & & Q_2 \; \geq \; 0, \; Q_2 + \bar Q_2  \; \geq \; 0,   \;\;\;  R_2 \; \geq  \delta I_m, 
%\right.
\enq
for some $\delta$ $>$ $0$,   the matrix Riccati equations \reff{Riccatilambda1}-\reff{Riccatigamma1}  admit unique solutions  $(\Lambda,\Gamma)$ $\in$  
$\Cc^1([0,T];\S_+^d)$ $\times$ $\Cc^1([0,T];\S_+^d)$,  and then $U_t,V_t$ in \reff{defUVSZ} are symmetric positive definite matrices, i.e. 
lie in $\S_{>+}^m$ for all $t$ $\in$ $[0,T]$.  The expression in \reff{opticontrolLQ} of the optimal control  extends then to the case of stochastic LQ McKean-Vlasov control problem  the 
feedback form obtained in \cite{yon13} for LQ McKean-Vlasov without  common noise, i.e. $\sigma_0$ $=$ $0$. 
%We shall see in the following  application  an example  with  explicit solutions where condition \reff{condmatrice} is not satisfied.  
}
\ep
\end{Remark}

\vspace{2mm}

\noindent {\bf Example: Interbank systemic risk model} 

\vspace{1mm}

\noindent We consider a model of inter-bank borrowing and lending studied in \cite{caretal14} where the log-monetary reserve of each bank in the asymptotics when the number of banks tend to infinity, is governed by 
the McKean-Vlasov equation:
\beq
dX_t &=&  \big[ \kappa(\E[X_t | W^0] - X_t)  + \alpha_t(X_t) ] dt  \nonumber \\
& & \;\;\; + \;  (\sigma_0 + \sigma_1 X_t)(\sqrt{1-\rho^2} dB_t + \rho dW_t^0), \; X_0 \; = \;  x_0  \in   \R.     \label{logX}
\enq 
Here, $\kappa$ $\geq$ $0$ is the rate of mean-reversion in the interaction from borrowing and lending between the banks,  $\sigma_0$ $>$ $0$, $\sigma_1$ $\in$ $\R$ are the affine coefficients of the volatility of the bank reserve,  
and there is a common noise $W^0$ for all the banks.  This is a slight extension of the model considered in \cite{caretal14} where $\sigma_1$ $=$ $0$. 
Moreover, all  banks  can  control  their  rate of borrowing/lending to a central bank with the same feedback policy $\alpha$  in order to minimize  a cost functional of the form 
\beqs
J(\alpha) &=&  \E \Big[ \int_0^T \Big( \frac{1}{2} \alpha_t(X_t)^2  - q \alpha_t(X_t) ( \E[X_t|W^0]-X_t)+ \frac{\eta}{2} (\E[X_t| W^0] - X_t)^2 \Big) dt   \\
& & \;\;\;\;\; + \;    \frac{c}{2}  (\E[X_T | W^0] - X_T)^2 \Big],   
\enqs
where $q$ $>$ $0$ is a positive parameter for the incentive to borrowing ($\alpha_t$ $>$ $0$) or lending ($\alpha_t$ $<$ $0$), and 
$\eta$ $>$ $0$, $c$ $>$ $0$ are positive parameters for penalizing departure from the average.  After square completion, we can rewrite the cost functional as
\beqs
J(\alpha) &=&  \E \Big[ \int_0^T \Big( \frac{1}{2} \tilde \alpha_t(X_t)^2   + \frac{\eta-q^2}{2} (\E[X_t| W^0] - X_t)^2 \Big) dt    +    \frac{c}{2}  (\E[X_T | W^0] - X_T)^2 \Big],   
\enqs
with $\tilde\alpha_t(X_t)$ $=$ $\alpha_t(X_t) - q(\E[X_t|W^0]-X_t)$.  This model fits into the framework of \reff{bsigLQ}-\reff{fgLQ} of the LQ stochastic  McKean-Vlasov problem with
\beqs
b_0=0, \; B = - (\kappa+q),\; \bar B= \kappa+q,\; C  =1, \\
D = \sigma_1\sqrt{1-\rho^2},  D_0 = \sigma_1 \rho, \; \bar D=F=\bar D^0=F^0 =0,  \; \vartheta=\sigma_0\sqrt{1-\rho^2}, \; \vartheta_0=\sigma_0 \rho,  \\
Q_2 =\frac{\eta-q^2}{2},\; \bar Q_2 =-\frac{\eta-q^2}{2},\; R_2 =\frac{1}{2}, 
%\; M_2  =-\frac{q}{2}, \; 
P_2=\frac{c}{2},\; \bar P_2=-\frac{c}{2}.
%q_1=\bar q_1=r_1 =p_1=\bar p_1=0.
\enqs
The Riccati system \reff{Riccatilambda1}-\reff{Riccatigamma1}-\reff{lingamma1}-\reff{linchi1}
for $(\Lambda(t), \Gamma(t), \gamma(t), \chi(t))$ is written in this case as
\begin{equation} \label{ODEex}
\left\{
\begin{array}{rclrccc}
\Lambda'(t)-2(\kappa +q - \frac{\sigma_1^2}{2}) \Lambda(t) - 2 \Lambda^2(t)  +   \frac{1}{2} (\eta-q^2) &=&0, & & \Lambda(T) &=& \frac{c}{2},\\
\Gamma'(t) - 2 \Gamma^2(t) + \sigma_1^2\rho^2 \Gamma(t) + \sigma_1^2(1-\rho^2)\Lambda(t)  &=&0, & & \Gamma(T) &=& 0,\\
\gamma'(t) -2 \Gamma(t)\gamma(t)  + 2 \sigma_0\sigma_1 \rho^2 \Gamma(t) + 2 \sigma_0\sigma_1 (1-\rho^2) \Lambda(t)  &=&0, & & \gamma(T) &= &0,\\
\chi'(t)-\frac{1}{2} \gamma^2(t) + \sigma_0^2 \rho^2 \Gamma(t) +\sigma_0^2 (1-\rho^2) \Lambda(t) &=&0, & &\chi(T)& =& 0.
\end{array}
\right.
\end{equation}
Assuming that $q^2$ $\leq$ $\eta$, the explicit solution to the Riccati equation for $\Lambda$ is given by
\beqs
%\chi(t) &=& \sigma^2(1-\rho^2) \int_t^T \Lambda(s) ds,  \\
\Lambda(t) &=& \frac{1}{2} \frac{ (\eta-q^2) \big(  e^{(\delta^+-\delta^-)(T-t)}-1\big)  +  c\big( \delta^+e^{(\delta^+-\delta^-)(T-t)} - \delta^- \big)}
{  c \big(e^{(\delta^+-\delta^-)(T-t)}-1\big)  + \delta^+  -  \delta^- e^{(\delta^+-\delta^-)(T-t)}  } \;\; > \; 0, 
\enqs 
where we set
\beqs
\delta^\pm &=& - \big(\kappa +q  - \frac{\sigma_1^2}{2}\big) \pm \sqrt{ \big(\kappa +q  - \frac{\sigma_1^2}{2}\big)^2 + \eta-q^2}. 
\enqs 
Since $\Lambda$ $\geq$ $0$, there exists a unique solution to the Riccati equation for $\Gamma$, and then $\gamma$, and finally $\chi$ are determined  the linear ordinary differential equations in \reff{ODEex}. 
Moreover,  the functions $(U_t,V_t,Z_t,Y_t)$ in \reff{defUVSZ} are explicitly given by:   $U_t$ $=$ $V_t$ $=$ $\frac{1}{2}$ (hence $>$ $0$),  
$S_t$ $=$ $\Lambda(t)+\frac{q}{2}$, $Z_t$ $=$ $\Gamma(t)$, $Y_t$ $=$ $\gamma(t)$. Therefore,  the optimal  control is given in  feedback form from  \reff{opticontrolLQ} by 
\beq 
\alpha_t^*(X_t^*) & = & a^*(t,X_t^*, \P_{_{X_t^*}})  \nonumber \\
&=& - (2\Lambda(t)+q)(X_t^* -\E[X_t^*| W^0]) - 2 \Gamma(t) \E[X_t^*| W^0] - \gamma(t), \label{controlsys}
\enq
where $X^*$ is the optimal  log-monetary reserve controlled by the rate  of borrowing/lending $\alpha^*$.  Moreover, denoting by $\bar X_t^*$ $=$ $\E[X_t^*|W^0]$ the conditional mean of the optimal log mo\-netary reserve, we see  that 
$\E[\alpha_t^*(X_t^*)|W^0]$ $=$ $-2\Gamma(t)\bar X_t^*-\gamma(t)$, and thus $\bar X^*$ is given from \reff{logX} by
\beqs
d\bar X_t^* &=& - \big(2 \Gamma(t) \bar X_t^* + \gamma(t) \big) dt + (\sigma_1 \bar X_t^*  + \sigma_0) \rho dW_t^0, \;\;\; \bar X_0^* = x_0. 
\enqs
When $\sigma_1$ $=$ $0$, we have $\Gamma(t)$ $=$ $\gamma(t)$ $=$ $0$, hence $\bar X_t^*$ $=$ $x_0$ $+$ $\sigma_0 \rho W_t^0$, 
and we  retrieve the expression found in \cite{caretal14} by sending the number of banks $N$ to infinity in their 
formula for the optimal control of the borro\-wing/lending rate:  
\beqs
\alpha_t^*(X_t^*)  &=& - (2\Lambda(t)+q)(X_t^* - x_0 - \sigma\rho W_t^0), \;\;\; 0 \leq t \leq T. 
\enqs

\end{document}